\documentclass[12pt]{article}

\usepackage{amstext    }
\usepackage{amsthm    }
\usepackage{a4}
\usepackage[mathscr]{eucal}
\usepackage{mathrsfs}

\usepackage{amsmath}
\usepackage{amssymb}
\usepackage{amscd}

\newtheorem{theorem}{Theorem}[section]
\newtheorem{definition}[theorem]{Definition}
\newtheorem{proposition}[theorem]{Proposition}
\newtheorem{corollary}[theorem]{Corollary}
\newtheorem{lemma}[theorem]{Lemma}
\newtheorem{remark}[theorem]{Remark}

\newtheorem{example}[theorem]{Example}
\newtheorem{examples}[theorem]{Examples}

\newcommand{\cali}[1]{\mathscr{#1}}

\newcommand{\volume}{{\rm vol}}

\newcommand{\supp}{{\rm supp}}

\newcommand{\dist}{{\rm dist}}

\newcommand{\loc}{{loc}}
\newcommand{\ddc}{{dd^c}}

\newcommand{\ddbar}{{\partial\overline\partial}}

\newcommand{\Gr}{{\rm Gr}}

\newcommand{\Jac}{{\rm Jac}}

\renewcommand{\Re}{{\rm Re}}
\renewcommand{\Im}{{\rm Im}}

\newcommand{\Cc}{\cali{C}}

\newcommand{\Fc}{\cali{F}}

\newcommand{\Lc}{\cali{L}}

\newcommand{\Pc}{\cali{P}}

\newcommand{\Lam}{{\rm Lam}}

\newcommand{\FS}{{\rm FS}}

\newcommand{\C}{\mathbb{C}}

\newcommand{\Z}{\mathbb{Z}}
\newcommand{\R}{\mathbb{R}}
\newcommand{\T}{\mathbb{T}}
\newcommand{\B}{\mathbb{B}}
\newcommand{\G}{\mathbb{G}}
\newcommand{\U}{\mathbb{U}}

\newcommand{\W}{\mathbb{W}}
\renewcommand{\S}{\mathbb{S}}
\renewcommand{\P}{\mathbb{P}}

\newcommand{\GG}{\mathbf{G}}

\renewcommand{\L}{\mathcal{L}}

\title{Equidistribution of saddle periodic points for H\'enon-type automorphisms of $\C^k$}

\author{Tien-Cuong Dinh and Nessim Sibony}

\begin{document}

\maketitle

\begin{abstract}
In this paper, we prove the equidistribution of saddle periodic points for H\'enon-type automorphisms of $\C^k$ with respect to it equilibrium measure. A general strategy to obtain equidistribution properties in any dimension is presented.  It is based on our recent 
 theory of densities for positive closed currents.  Several fine properties of dynamical currents are also proved.
\end{abstract}

\noindent
{\bf Classification AMS 2010:} 32U, 37F, 32H50.


\section{Introduction} \label{intro}

Let $f$ be a polynomial automorphism of $\C^k$. We extend it to $\P^k$ as a bi-rational self-map that we still denote by $f$. 
Let $I_+$ denote the indeterminacy set of $f$ and $I_-$ the one for the inverse $f^{-1}$ of $f$. They are contained in the hyperplane at infinity $H_\infty:=\P^k\setminus\C^k$. We assume that $f$ is not an automorphism of $\P^k$ because otherwise its dynamics is elementary. So the indeterminacy sets $I_+$ and $I_-$ are non-empty. The following notion was introduced by the second author in \cite{Sibony}.

\begin{definition} \rm
We say that $f$ is a {\it H\'enon map} or a {\it regular automorphism} if 
$$I_+\cap I_-=\varnothing.$$ 
\end{definition}

The interesting point here is that the last condition is quite simple to check and H\'enon automorphisms form a rich family of non-uniformly hyperbolic dynamical systems for which we can develop a satisfactory theory. In dimension 2, all dynamically interesting automorphisms of $\C^2$ are conjugated to H\'enon maps, see Friedland-Milnor \cite{FriedlandMilnor} and also \cite{DS15, Fornaess, FornaessWu}. 
Assume now that $f$ is a H\'enon map on $\C^k$.
We first recall some basic properties of  $f$ and refer to the papers by de Th\'elin \cite{deThelin} and the authors \cite{Dinh,DS10,Sibony} for details.

Let $d_\pm$ denote the algebraic degrees of $f^{\pm 1}$. There is an integer $1\leq p\leq k-1$ such that $\dim I_+=k-p-1$, $\dim I_-=p-1$ and $d_+^p=d_-^{k-p}\geq 2$. We define the {\it main dynamical degree} of $f$ as $d:=d_+^p=d_-^{k-p}$. This is also the main dynamical degree of $f^{-1}$ and the topological entropies of $f$ and $f^{-1}$ are both equal to $\log d$. 
The restrictions of $f$ and its inverse to the hyperplane at infinity $H_\infty$ satisfy  
$$f(H_\infty\setminus I_+)=f(I_-)=I_- \quad \text{and}\quad  f^{-1}(H_\infty\setminus I_-)=f^{-1}(I_+)=I_+.$$ 

Let $K_+$ (resp. $K_-$)  be the set of points $z\in\C^k$ such that the orbits $(f^n(z))_{n\geq 0}$ (resp. $(f^{-n}(z))_{n\geq 0}$) are bounded in $\C^k$. They are closed in $\C^k$ and we have 
$$\overline K_\pm=K_\pm\cup I_\pm.$$ 
The indeterminacy sets $I_-$ and $I_+$ are attracting respectively for $f$ and $f^{-1}$. Their basins are equal to $\P^k\setminus \overline K_+$ and $\P^k\setminus \overline K_-$. The intersection $K:=K_+\cap K_-$ is compact in $\C^k$. It is the set of points $z$ whose entire orbits $(f^n(z))_{n\in \Z}$ are bounded in $\C^k$. 

Let $\omega_\FS$ denote the Fubini-Study $(1,1)$-form on $\P^k$ so normalized  that the integral of the top power $\omega_\FS^k$ on $\P^k$ is 1. The following weak limits  exist
$$\tau_\pm:=\lim_{n\to\infty} d_\pm^{-n} (f^{\pm n})^*(\omega_\FS)$$
and define two positive\footnote{Throughout the paper, the positivity of $(p,p)$-currents is in the strong sense.} closed $(1,1)$-currents of mass 1 on $\P^k$. We have $f^*(\tau_+)=d_+\tau_+$ and $f_*(\tau_-)=d_-\tau_-$. The currents $\tau_+$ and $\tau_-$ have H\"older continuous local potentials outside $I_+$ and $I_-$  respectively. They are called the {\it Green $(1,1)$-currents} of $f$ and $f^{-1}$. 

The positive closed currents $T_+:=\tau_+^p$ and $T_-:=\tau_-^{k-p}$ are respectively the main Green currents of $f$ and $f^{-1}$. The current $T_+$ is the unique positive closed $(p,p)$-current of mass 1 in $\P^k$ with support in $\overline K_+$ and  
the current $T_-$ is the unique positive closed $(k-p,k-p)$-current of mass 1 in $\P^k$ with support in $\overline K_-$.

The wedge-product
$$\mu:=T_+\wedge T_-=\tau_+^p\wedge\tau_-^{k-p}$$
is a well-defined invariant probability measure with support in $K$. It turns out that $\mu$ is the unique invariant probability measure of maximal entropy $\log d$ of $f$ and $f^{-1}$. The measure $\mu$ is moreover exponentially mixing and hyperbolic. It is called the {\it Green measure or equilibrium measure} of $f$ and $f^{-1}$. 

In this paper, we give the proof that saddle periodic points are equidistributed with respect to $\mu$. Denote by $P_n$ the set of periodic points of period $n$ of $f$ in $\C^k$ and $SP_n$ the set of those which are saddles. For any number $0<\epsilon<1$, denote by $SP_n^\epsilon$ the set of saddle periodic points $a$ of period $n$ in $\C^k$ such that the differential $Df^n(a)$ admits exactly $p$ eigenvalues of modulus larger than $(d_+-\epsilon)^{n/2}$ and    $k-p$ eigenvalues of modulus smaller than $(d_--\epsilon)^{-n/2}$. Here the eigenvalues are counted with multiplicity. They do not depend on the choice of coordinate system on $\C^k$.
We have the following theorem, see Bedford-Lyubich-Smillie \cite{BLS2} for the case of dimension $k=2$ and \cite{Lee} for a $p$-adic version independently obtained by Lee. The main result by Lee may offer an arithmetic approach this this problem by taking $p\to\infty$. 

\begin{theorem} \label{th_main}
Let $f,d, \mu,P_n,SP_n$ and $SP_n^\epsilon$ be as above. Let $Q_n$ be either $P_n, SP_n$ or $SP_n^\epsilon$. Then 
$$ d^{-n} \sum_{a\in Q_n} \delta_a\to \mu$$
as $n$ goes to infinity, where $\delta_a$ denotes the Dirac mass at $a$.
\end{theorem}

The proof of this result is developed in the rest of the paper. A key point is the use of our theory of densities of positive closed currents developed in \cite{DS12}. We refer to that paper for basic notations and results concerning tangent cones, the notion of density and the intersection of currents in a weak sense. We will describe below our strategy which, as far as we know, is the first approach to obtain the equidistribution of periodic points for a non-uniformly hyperbolic holomorphic system with arbitrary numbers of stable and unstable directions. The main ideas are quite general and can be adapted to other meromorphic dynamical systems and other questions, see also Remark \ref{rk_final} below.


Let $\Delta$ denote the diagonal of $\P^k\times\P^k$ and $\Gamma_n$ denote the compactification of the graph of $f^n$ in $\P^k\times\P^k$. The set $P_n$ can be identified with the intersection of $\Gamma_n$ and $\Delta$ in $\C^k\times\C^k$. The dynamical system associated with the map $F:=(f,f^{-1})$ on $\P^k\times\P^k$ is similar to the one associated with H\'enon-type maps on $\P^k$. It was used by the first author in \cite{Dinh} in order to obtain the exponential mixing of $\mu$ on $\C^k$. 
Observe that $\Gamma_n$ is the pull-back of $\Delta$ or $\Gamma_1$ by $F^{n/2}$ or $F^{(n-1)/2}$. 
So a property similar to the uniqueness of the main Green currents mentioned above allows us to prove that the positive closed $(k,k)$-current $d^{-n}[\Gamma_n]$ converges to the main Green current of $F$ which is equal to
$T_+\otimes T_-$. Therefore, since the measure $\mu=T_+\wedge T_-$ can be identified with $[\Delta]\wedge (T_+\otimes T_-)$,  Theorem \ref{th_main} is equivalent to
$$\lim_{n\to\infty} [\Delta]\wedge d^{-n}[\Gamma_n]=
[\Delta]\wedge \lim_{n\to\infty}  d^{-n}[\Gamma_n]$$
on $\C^k\times\C^k$. 
So our result requires the development of a good intersection theory in any dimension.

The typical difficulty is illustrated in the following example. Consider $\Delta'$ the unit disc in $\C\times \{0\}\subset \C^2$ and $\Gamma_n'$ the graph of the function $x\mapsto x^{d^n}$ over $\Delta'$. The currents $d^{-n}[\Gamma_n']$ converge to a current on the vertical boundary of the unit bidisc in $\C^2$ while their intersection with $[\Delta']$ is the Dirac mass at 0. So we have
$$\lim_{n\to\infty} [\Delta']\wedge d^{-n}[\Gamma_n']\not =
[\Delta']\wedge \lim_{n\to\infty}  d^{-n}[\Gamma_n'].$$
We see in this example that $\Gamma_n'$ is tangent to $\Delta'$ at 0 with maximal order. We can perturb $\Gamma_n'$ in order to get manifolds which intersect $\Delta'$ transversally but the limit of their intersections with $\Delta'$ is still equal to the Dirac mass at 0.  In fact, this phenomenon is due to the property that some tangent lines to $\Gamma_n'$ are too close to tangent lines to $\Delta'$. 

It is not difficult to construct a map $f$ such that $\Gamma_n$ is tangent or almost tangent to $\Delta$ at some points for every $n$. In order to handle the main difficulty in our problem, the strategy is to show that the almost tangencies become negligible when $n$ tends to infinity. This property is translated in our study into the fact that a suitable density for positive closed currents vanishes. Then, a geometric approach developed in \cite{Dinh2} allows us to obtain the result. We will give the details in the second part of this article. We explain now the notion of density of currents in the dynamical setting and then develop the theory in the general setting of arbitrary positive closed currents. 

Let $\Gr(\P^k\times\P^k,k)$ denote the Grassmannian bundle over $\P^k\times\P^k$ where each point corresponds to a pair $(x,[v])$ of a point $x\in\P^k\times\P^k$ and the direction $[v]$ 
of a simple tangent $k$-vector $v$ of $\P^k\times\P^k$ at $x$.   Let $\widetilde\Gamma_n$ denote  the set of points $(x,[v])$ in $\Gr(\P^k\times\P^k,k)$
with $x\in\Gamma_n$ and $v$ a $k$-vector not transverse to $\Gamma_n$ at $x$. Let $\widehat\Delta$ denote the lift of $\Delta$ to 
$\Gr(\P^k\times\P^k,k)$, i.e.
the set of points $(x,[v])$  with $x\in\Delta$ and $v$ tangent to $\Delta$.  The intersection $\widetilde\Gamma_n\cap \widehat\Delta$ corresponds to the non-transverse points of intersection between $\Gamma_n$ and $\Delta$. Note that $\dim \widetilde \Gamma_n+\dim \widehat\Delta$ is smaller than the dimension of $\Gr(\P^k\times\P^k,k)$ and the intersection of subvarieties of such dimensions are generically empty. Analogous construction can be done for the manifolds $\Gamma_n'$ and $\Delta'$ given above.

We show that the currents $d^{-n}[\widetilde\Gamma_n]$ cluster on some positive closed current $\widetilde \T_+$ on $\Gr(\P^k\times\P^k,k)$. It can be obtained from the current $\T_+:=T_+\otimes T_-$ on $\P^k\times\P^k$ in a similar way as in the construction of $\widetilde\Gamma_n$. One can also construct $\widetilde\T_+$  by lifting $\T_+$ to a positive closed current $\widehat\T_+$ of the same dimension in $\Gr(\P^k\times\P^k)$ and then transforming it into $\widetilde\T_+$ via some incidence manifold. 
Using a theorem due to de Th\'elin \cite{deThelin2} on the hyperbolicity of $\mu$ we show that the density between $\widetilde \T_+$ and $\widehat\Delta$ vanishes. This property says that almost tangencies are negligible when $n$ goes to infinity. The above example with $\Gamma_n'$ and $\Delta'$ is an illustration of the opposite situation.  

Let us be more precise. The vanishing of the density between $\widetilde\T_+$ and $\widehat\Delta$ implies that the mass of $\widetilde \Gamma_n$ in a small enough neighbourhood of $\widehat \Delta$ is smaller than $\epsilon d^n$ for any given small constant $0<\epsilon<1$ when $n$ is large enough. It follows that for some projection close to the projection $(x,y)\mapsto x-y$ from $\C^k\times\C^k$ to $\C^k$ the size of the ramification locus of $\Gamma_n$ over a small neighbourhood of $0\in\C^k$ is smaller than $\epsilon d^n$. On the other hand, with respect to this projection, $\Gamma_n$ is a ramified covering over $\C^k$ of degree approximatively $d^n$. An argument \`a la Hurwitz permits to construct almost $d^n$ graphs $\Gamma_n^{(j)}$ contained in $\Gamma_n$ over a small neighbourhood of $0\in\C^k$. 

Each graph $\Gamma_n^{(j)}$ intersects $\Delta$ at a unique point corresponding to a periodic point of $f$. If a sequence of such graphs converges in the sense of currents, it converges uniformly. Therefore, we control the convergence of a large part of the intersection $\Gamma_n\cap \Delta$ thanks to the convergence of $\Gamma_n$. This together with some standard arguments imply the identity
$$\lim_{n\to\infty} [\Delta]\wedge d^{-n}[\Gamma_n]=
[\Delta]\wedge \lim_{n\to\infty}  d^{-n}[\Gamma_n]$$
which is equivalent to Theorem \ref{th_main}.

The dynamical setting enters into the picture, first because in the above arguments we use that the graphs $\Gamma_n$ are horizontal with respect to the projection $(x,y)\mapsto x-y$. Near the diagonal $\Delta$, they are contained in a fixed box along $\Delta$. The other more serious point concerns the delicate computation of the density between $\widetilde\T_+$ and $\widehat\Delta$.  
Roughly speaking, we want to show that $\widetilde\T_+$ is not concentrated near $\widehat\Delta$ as we can observe positive closed currents near a point with positive Lelong number.

For $\mu$-almost every point $z\in\C^k$, denote by $E_s(z)$ and $E_u(z)$ the stable and unstable tangent subspaces for $f$ and $\mu$  at $x$. Since $\mu$ is hyperbolic with $p$ positive Lyapounov exponents and $k-p$ negative ones, we have $\dim E_s(z)=k-p$, $\dim E_u(z)=p$ and $E_s(z)\cap E_u(z)=\varnothing$. Denote by $\Pi$ the canonical projection from $\Gr(\P^k\times\P^k,k)$ to $\P^k\times\P^k$. Let $\GG(z)$ be the set of points $(z,z,[v])$ in $\Pi^{-1}(z,z)$ such that $v$ is not transverse to the vector space 
$E_s(z)\times E_u(z)$. Since the last vector space is transverse to $\Delta$, the varieties $\GG(z)$ and $\widehat\Delta$ are disjoint. 

We can show that the intersection of $\widetilde\T_+$ and $\Pi^{-1}(\Delta)$ is, in a weak sense of currents defined later,  equal to the average,  with respect to $\mu$, of the currents of integration on $\GG(z)$. This delicate property is basically due to the hyperbolicity of $\mu$ and Oseledec's  theorem. We see roughly that $\widetilde\T$ crosses $\Pi^{-1}(\Delta)$ through the varieties $\GG(z)$ which are disjoint from $\widehat\Delta$. This is the key point to get the vanishing of the mentioned density between $\widetilde\T_+$ and $\widehat\Delta$. The proof requires also several geometric properties of $\widetilde\T_+$ and $\widehat \T_+$. They are of independent interest.

\medskip
\noindent
{\bf Acknowledgements. } The paper was partially written during the visit of the
first author at the Shanghai Center for Mathematical Sciences. He would like to thank the institute, Yi-Jun Yao and Weiping Zhang for 
their great hospitality.

\section{Preliminary on currents and cohomology} \label{section_coh}

In this section, we will give  some properties of positive closed currents and some properties of the action of meromorphic maps on currents and on Hodge cohomology. They will be used to overcome technical difficulties in the proof of our main result. 
Some of them are of independent interest. Though the strategy should work in a quite general setting, it is already involved for H\'enon maps on $\C^k$. So we will limit ourself to the simplest situation required for H\'enon maps  and refer to  \cite{DNT, DS3} for the case of general meromorphic maps on compact K\"ahler manifolds.

\medskip\noindent
{\bf Slicing theory for currents.} 
We will discuss slicing theory in the case of positive closed currents. It will be applied later to woven positive closed currents. Recall that Federer's slicing theory can be applied to positive closed currents and to currents defined by  complex varieties which are not necessarily closed, see Federer \cite{Federer}. The later case allows to extend the slicing theory to woven currents, not necessarily closed, that we will introduce in the next section. 

Let $\pi:V\to W$ be a proper surjective holomorphic map between K\"ahler manifolds of dimension $l$ and $m$ respectively with $l\geq m$. Let $T$ be a positive closed $(p,p)$-current on $V$ of bi-dimension at least equal to $(m,m)$, i.e. $p\leq l-m$. Then for almost every $y\in W$, the slice $\langle T|\pi|y\rangle$ exists and is a positive closed current of bi-dimension $(l-p-m,l-p-m)$ and of bi-degree $(p+m,p+m)$ with support in $\pi^{-1}(y)$.  This current is obtained as the limit of a sequence of currents of the form $T\wedge \pi^*(\phi_n)$, where $(\phi_n)$ is a suitable sequence of positive $(m,m)$-forms on $W$ which converges weakly to the Dirac mass at $y$. Moreover, if $\phi$ is a continuous $(m,m)$-form on $W$ and $\psi$ is a continuous $(l-p-m,l-p-m)$-form with compact support on $V$ then
$$\langle T,\pi^*(\phi)\wedge\psi\rangle = \int_{y\in W} \langle T|\pi|y\rangle(\psi) \phi(y).$$

Assume now that $\pi$ is proper on a closed subset $G$ of $V$ which contains the support of $T$. We developed in \cite{DS11} a simpler slicing theory with some advantages that we recall below. Let $\psi$ be an $(l-p-m,l-p-m)$-form of class $\Cc^2$ on $V$. Assume it is a real form; otherwise, we can consider separately its real and imaginary parts. Then $\langle T|\pi|y\rangle(\psi)$ is equal almost everywhere on $W$ to a {\it d.s.h.} function, i.e. locally a difference of two psh functions, see Remark 2.2.6 in \cite{DS11}. If the mass of $T$ is bounded and $\psi$ is fixed, we can take these local psh functions in a suitable compact family of psh functions. Here is a consequence of this property.

\begin{proposition} \label{prop_slice_C0}
Let $V, W,\pi$ and $G$ be as above. Let $T_n$ be a sequence of positive closed $(p,p)$-currents on $V$ converging to a positive closed $(p,p)$-current $T$ with $p\leq l-m$.  Assume that all $T_n$ are supported by $G$. Then there is a subsequence $T_{n_i}$ such that for almost every $y\in W$ we have
$$\lim_{i\to\infty} \langle T_{n_i}|\pi|y\rangle =\langle T|\pi|y\rangle.$$ 
\end{proposition}
\proof
Recall that the map embedding psh functions in $L^1_\loc$ is compact. 
From the d.s.h. property we just mentioned, for every $\Cc^2$ test function $\psi$, the sequence  of d.s.h. functions $\langle T_{n}|\pi|y\rangle(\psi)$ on $W$ converges to $\langle T|\pi|y\rangle(\psi)$ in $L^1_\loc$. 
Recall also that any convergent sequence in $L^1_\loc$ admits a subsequence which converges almost everywhere. We deduce from the property of slicing mentioned above that  for a given $\psi$ in $\Cc^2$, there is a subsequence $T_{n_i}$ such that we have almost everywhere on $W$
$$\lim_{i\to\infty} \langle T_{n_i}|\pi|y\rangle(\psi) =\langle T|\pi|y\rangle(\psi).$$ 

Using the standard diagonal process, we can have such a property of a suitable sequence $(n_i)$  for any fixed countable family of test $\Cc^2$ forms $\psi$. Choose such a family which is dense in the space $\Cc^0$. The last identity applied to $\psi$ in this family implies the identity in the proposition. 
\endproof

Our slicing theory does not apply to non-closed currents. Consider the case where $T$ is the current of integration on a complex variety $\Lambda$ of dimension $l-p$, immersed in $V$, whose $2(l-p)$-dimensional volume is finite in each compact subset of $V$. Then for almost every $y\in W$  the intersection $\Lambda\cap \pi^{-1}(y)$ is either empty or a variety of dimension $l-p-m$ with finite $2(l-p-m)$-dimensional volume in each compact subset of $V$ and the intersection is transverse at almost every point.  
Moreover, the slice $\langle T|\pi|y\rangle$ in the classical sense of Federer and is equal to 
the current of integration on $\Lambda\cap \pi^{-1}(y)$. 
For simplicity, we only consider generic $y$ satisfying these properties. 
 Note that in the case where the restriction of $\pi$ to $V$ has rank strictly smaller than $m$, then for almost every $y$ we have $\langle T|\pi|y\rangle=0$.

We describe now a situation which will be used later. 
Consider the case where $V=V'\times \P^k$ with $V'$ compact K\"ahler of dimension $l-k$. If $H$ is a generic projective subspace of dimension $r$ in $\P^k$, we will define the intersection of $T$ with $V'\times H$, in the slicing theory sense. Let $\G$ denote the Grassmannian which parametrizes the family of such $H$. For $\xi\in\G$, denote by $H_\xi$ the corresponding projective subspace in $\P^k$. Let $\Sigma$ be the incidence manifold defined by
$$\Sigma:=\{(x,\xi) \text{ with } x\in H_\xi\subset\P^k \text{ and } \xi\in\G\}\subset \P^k\times\G.$$
Let $\pi_1,\pi_2$ denote the natural projections from $V'\times \Sigma$ onto $V=V'\times\P^k$ and $\G$. They are submersions. So the positive closed current $\pi_1^*(T)$ is well-defined.
Define  the intersection $T\wedge [V'\times H_\xi]$ as the push-forward of the slice $\langle \pi_1^*(T)|\pi_2| \xi\rangle$ by $\pi_1$ to $V$. This slice exists for generic $\xi\in \G$. We will be concerned with the case where $V'$ is a product of $\P^k$ and a Grassmannian bundle over $\P^k$.

\medskip\noindent
{\bf Some general norm controls on pull-back operators.}
The manifolds involved in our study are bi-rational to products of projective spaces and Grassmannians. We will sometimes reduce the problem to the case of such products where the structure of Hodge cohomology groups is simpler.
The following proposition allows us to do this reduction.  It is a direct consequence of the main theorem in \cite{DS3}. 

\begin{proposition} \label{prop_conj_map}
Let $\pi:X'\to X$ be a bi-meromorphic map between two compact K\"ahler manifolds of dimension $k$. Let $f,g:X\to X$ and $f':X'\to X'$ be dominant meromorphic maps. Assume that $f$ and $f'$ are conjugated: $f\circ \pi=\pi\circ f'$. 
Let $f^*$ and $g^*$ (resp. ${f'}^*$) denote the pull-back operators associated with $f$ and $g$ (resp. $f'$) acting on the Hodge cohomology group $H^{q,q}(X,\C)$ (resp. $H^{q,q}(X',\C)$)  for an integer  $0\leq q\leq k$.  
Then
there is a constant $A\geq 1$ depending only on $X,X',\pi$ and the norms on Hodge cohomologies of $X,X'$ such that 
$$A^{-1}\|{f'}^*\|\leq \|f^*\|\leq A\|{f'}^*\| \quad \text{and}\quad \|f^*\circ g^*\|\leq A\|f^*\|\|g^*\|.$$
\end{proposition}

Note that on a compact K\"ahler manifold with a fixed K\"ahler metric the mass of a positive closed current depends only on its cohomology class. So the last proposition is equivalent to a mass control for these currents under the action of meromorphic maps. We can apply it to $f^n, g^n$ and ${f'}^n$ instead of $f,g$ and $f'$ keeping the same constant $A$. 

We will also need the following proposition in order to control the mass of a current, see \cite{DS3} for a more general result.

\begin{proposition} \label{prop_bullet}
Let $f:X\to X$ be a dominant meromorphic map on a compact K\"ahler manifold of dimension $k$. Let $X_0$ be a non-empty Zariski open subset  of $X$ such that $f$ defines a bi-holomorphic map between $X_0$ and its image. Let $T$ be a positive closed $(q,q)$-current of mass $1$ on $X$ with  $0\leq q\leq k$. Then the pull-back $(f_{|X_0})^*(T)$ of $T$ to $X_0$ has finite mass. 
Denote by $f^\bullet(T)$ its extension by $0$ to $X$. Then $f^\bullet(T)$ is a positive closed $(q,q)$-current. Moreover, its mass is bounded by the norm of $f^*$ on $H^{q,q}(X,\C)$ times a constant $A$ which only depends on the K\"ahler metric  on $X$ and the norm on $H^{q,q}(X,\C)$. 
\end{proposition}

Note that $f^\bullet(T)$  depends on the choice of $X_0$ but in applications we often have a standard choice. We can apply the proposition to $f^n$ instead of $f$ with the same constant $A$. This constant does not depend on the choice of $X_0$. Note also that we have $(f^n)^\bullet=(f^\bullet)^n$ when $f$ defines an automorphism of $X_0$. 

\medskip
\noindent
{\bf Tensor products of  currents.} We will need criteria to check that a positive closed current is the tensor product of two other ones. 
Let $X$ and $Y$ be compact K\"ahler manifolds of dimension $k$ and $l$. For simplicity, assume that $H^{r,s}(X,\C)=0$ when $r\not=s$.
Let $\pi_X,\pi_Y$ be the canonical projections from $X\times Y$ onto its factors $X$ and $Y$ respectively. 
Fix K\"ahler forms $\omega_X,\omega_Y$  on $X$ and $Y$.
We have the following result.

\begin{proposition} \label{prop_current_tensor}
Let $T$ be a positive closed current on $X\times Y$ of bi-dimension $(s,s)$. Assume there is an integer $0\leq r\leq \min(k,s)$ such that $T\wedge \pi_X^*(\omega_X^{r+1})=0$ and $T\wedge \pi_Y^*(\omega_Y^{s-r+1})=0$. Assume also that  $R:=(\pi_X)_*\big(T\wedge \pi_Y^*(\omega_Y^{s-r})\big)$ is an extremal positive closed current of bi-dimension $(r,r)$ on $X$. Then there is a positive closed current $S$ of bi-dimension $(s-r,s-r)$ on $Y$ such that $T=R\otimes S$. 
\end{proposition}
\proof
Let $x$ and $y$ denote coordinates on $X$ and $Y$. Observe that each smooth $(s,s)$-form on $X\times Y$ can be written as a finite linear combination of forms of the following type or their conjugates:
$$\Phi:=h(x,y)\alpha(x)\wedge\beta(y)\wedge\Omega(x)\wedge \Theta(y),$$
where $h$ is a smooth positive function on $X\times Y$, $\alpha$ is a $(u,0)$-form on $X$, $\beta$ is a $(0,u)$-form on $Y$, $\Omega$ is a smooth positive $(v,v)$-form on $X$ and $\Theta$ is a smooth positive form of bi-degree $(s-u-v,s-u-v)$ on $Y$. 
We first prove the following claim. 

\medskip\noindent
{\bf Claim. } We have $\langle T,\Phi\rangle=0$ when $u\not=0$ or $v\not=r$.

\medskip

Observe that since $T$ is positive and $\omega_X$ is strictly positive, the hypotheses on $T$ imply that for any $(r+1,r+1)$-form $\phi$ on $X$, we have $T\wedge \pi_X^*(\phi)=0$. We then deduce that the same identity also holds  for every $(u,u)$-form $\phi$ with $u\geq r+1$.
Similarly, we have $T\wedge \pi_Y^*(\psi)=0$ for any $(v,v)$-form $\psi$ on $Y$ with $v\geq s-r+1$.

Observe also that 
for any constant $\lambda > 0$, the forms
$$i^{u^2} \big[\lambda\alpha(x)\wedge \overline{\alpha(x)}+\lambda^{-1}\overline{\beta(y)}\wedge\beta(y)\pm 2\Re \big(\alpha(x)\wedge\beta(y)\big)\big]$$
and
$$i^{u^2} \big[\lambda\alpha(x)\wedge \overline{\alpha(x)}+\lambda^{-1}\overline{\beta(y)}\wedge\beta(y)\pm 2\Im \big(\alpha(x)\wedge\beta(y)\big)\big]$$
are weakly positive. Since $T,\Omega,\Theta$ are positive and $h$ is bounded, we can bound $|\langle T,\Phi\rangle|$ by a constant times
$$\lambda \big\langle T, i^{u^2} \alpha(x)\wedge \overline{\alpha(x)}\wedge\Omega(x)\wedge \Theta(y)\big\rangle + \lambda^{-1} \big\langle T, i^{u^2} \overline{\beta(y)}\wedge\beta(y)\wedge\Omega(x)\wedge \Theta(y)\big\rangle.$$
For a good choice of $\lambda$, the last expression is equal to
$$ \big\langle T, i^{u^2} \alpha(x)\wedge \overline{\alpha(x)}\wedge\Omega(x)\wedge \Theta(y)\big\rangle^{1/2}\big\langle T, i^{u^2} \overline{\beta(y)}\wedge\beta(y)\wedge\Omega(x)\wedge \Theta(y)\big\rangle^{1/2}.$$
It is not difficult to see that the second factor vanishes when $v\not=r$ and the first factor vanishes when
$v=r$ and $u\not=0$. This completes the proof of the claim. 

\medskip

For a bi-degree reason, we easily deduce from the claim that $T$ vanishes on $\phi(x)\wedge \psi(y)$ for all smooth $(u,v)$-form $\phi$ on $X$ and $(s-u,s-v)$-form $\psi$ on $Y$ provided that $(u,v)\not=(r,r)$. Consider now a smooth $(r,r)$-form $\phi$ on $X$ and a smooth $(s-r,s-r)$-form $\psi$ on $Y$. 
Observe that $d(T\wedge \psi(y))=T\wedge d\psi(y)=0$ thanks to the last observation. 
So if $\psi$ is positive, then $T\wedge \psi(y)$ is a positive closed current. Its horizontal dimension with respect to the projection $\pi_Y$ is 0 in the sense that $T\wedge\psi(y)\wedge \omega_Y(y)=0$. 
Therefore, we can prove as in Lemma 3.3 in \cite{DS12} that $T\wedge\psi(y)$
can be disintegrated into positive closed currents on $X\times \{y\}$ with respect to a positive measure on $Y$. 
The push-forward of  $T\wedge \psi(y)$ to $X$ is bounded by a constant times $R$ since $\psi$ is bounded by a constant times $\omega_Y^{s-r}$.  Since $R$ is extremal, the above currents on $X\times \{y\}$ are proportional to $R$. The property holds without the positivity of $\psi$ since we can always write $\psi$ as a linear combination of positive forms. 

We deduce from the last property that 
 $$\langle T,\phi(x)\wedge\psi(y)\rangle = \langle R,\phi\rangle S(\psi),$$
 where $S$ is some continuous linear form, i.e. a current of bi-dimension $(s-r,s-r)$ on $Y$. The last identity also holds for 
 $(u,v)\not=(r,r)$ since in this case its both sides vanish. We deduce that $T=R\otimes S$ because the vector space generated by the  forms $\phi(x)\wedge\psi(y)$ is dense in the space of test $(s,s)$-forms on $X\times Y$. Since $T$ is positive closed, it is easy to check that $S$ is also positive and closed. This completes the proof of the proposition.
\endproof

In the dynamical setting, we will use the proposition below in order to check the hypotheses of the last result. 
Let $f:X\to X$ and $g:Y\to Y$ be bi-meromorphic maps. Assume that $X,Y$ are homogeneous and 
there are dense Zariski open sets $X_0\subset X$ and $Y_0\subset Y$ such that $f,g$ are automorphisms of $X_0$ and $Y_0$ respectively.  

Define the bi-meromorphic map $h:X\times Y\to X\times Y$ by $h(x,y):=(f(x),g(y))$. 
Assume there are constants $d>\delta>1$ and an integer $1\leq p\leq k-1$  such that $\|(f^n)^*\|=O(\delta^n)$ on $H^{q,q}(X,\C)$ 
for $q\not=p$ and $\|(f^n)^*\|=O(d^n)$ on $H^{p,p}(X,\C)$ as $n\to\infty$.
Fix also a constant $d'>1$ such that $\|(g^n)^*\|=O({d'}^n)$ on $H^{q,q}(Y,\C)$ for every $q$.

\begin{proposition} \label{prop_product_map}
Let $S_n$ be a sequence of positive closed currents of bi-dimension $(s,s)$ on $X\times Y$ with bounded mass. Let $T$ be a limit value of the sequence of currents $d^{-n}{d'}^{-n} (h^n)^\bullet(S_n)$. Then we have $T\wedge \pi_X^*(\omega_X^{r+1})=0$ and $T\wedge \pi_Y^*(\omega_Y^{s-r+1})=0$ for $r:=k-p$.
\end{proposition}

A $(q,q)$-class in the cohomology group of a compact K\"ahler manifold is said to be 
{\it pseudo-effective} if 
 it contains a positive closed current.  
 So we can define a partial order relation between real classes: we have $c\leq c'$ if $c'-c$ is pseudo-effective.
 We also say that a class is {\it strictly positive} if it is larger than or equal to the class of a strictly positive closed form. 
Recall also that on homogeneous manifolds, any positive closed currents can be approximated by smooth positive closed forms in the same cohomology class. These forms can be obtained using a convolution with holomorphic automorphisms close to the identity, see e.g. \cite{Hu}.

Recall also the K\"unneth formula in our case where $H^{r,s}(X,\C)=0$ for $r\not=s$, see \cite{Voisin}. We have the following canonical decomposition of the Hodge cohomology on $X\times Y$
$$H^{q,q}(X\times Y,\C)=\sum_r H^{r,r}(X,\C)\times H^{q-r,q-r}(Y,\C).$$
Here, for simplicity, we set $H^{r,r}(X,\C)=0$ if either $r<0$ or $r>\dim X$ and we apply the same convention to all manifolds. 

\medskip\noindent
{\bf Proof of Proposition \ref{prop_product_map}.} 
Since $(h^n)^*$ preserves the K\"unneth decomposition for $X\times Y$, it is not difficult to see that 
$\|(h^n)^*\|=O(d^n{d'}^n)$ on the Hodge cohomology of $X\times Y$. Note that $(h^n)^\bullet$ is not compatible with the action of $h^n$ on cohomology. 

Since $X\times Y$ is homogeneous, $S_n$ can be approximated by smooth positive forms $S_n^{(j)}$ in the same cohomology class. We deduce that $d^{-n}{d'}^{-n} (h^n)^\bullet(S_n)$ is smaller than or equal to all limit values of $d^{-n}{d'}^{-n} (h^n)^*(S_n^{(j)})$ when $j\to\infty$.
If $c_n$ denotes the cohomology class of $S_n$,  the class of   $d^{-n}{d'}^{-n} (h^n)^\bullet(S_n)$ is smaller than or equal to 
 $d^{-n}{d'}^{-n} (h^n)^*(c_n)$. So we only have to check that all limit values of $d^{-n}{d'}^{-n} (h^n)^*(c_n)$  belong to the component 
 $H^{p,p}(X,\C)\times H^{r+l-s,r+l-s}(Y,\C)$. 
 
 Since $(h^n)^*$ preserve the K\"unneth decomposition, it is enough to verify that the sequence of operators $d^{-n}{d'}^{-n} (h^n)^*$ converges to 0 on   $H^{q,q}(X,\C)\times H^{q',q'}(Y,\C)$ for $q\not=p$ and for every $q'$. This is clear because $(h^n)^*$ is the product of the operator $(f^n)^*$ acting on $H^{q,q}(X,\C)$ and the operator $(g^n)^*$ acting on $H^{q',q'}(Y,\C)$. The norm of first operator is equal to  $O(\delta^n)$ and the norm of the second one is $O({d'}^n)$.    The result follows.
 \hfill $\square$

\medskip\noindent
{\bf Some norm control on pull-back operators.}
Let $f:X\to X$ be a bi-meromorphic map on a compact K\"ahler manifold $X$ of dimension $k$. Let $X_0$ be a non-empty Zariski open subset of $X$ such that $f$ is a bi-holomorphism between $X_0$ and its image. Assume that 
$H^{r,s}(X,\C)=0$ for $r\not=s$.
Let $1\leq p\leq k-1$, $n_0\geq 1$ be integers and let $1\leq \delta<d$ be real numbers such that 
$\|(f^n)^*\|=O(\delta^n)$ on $H^{q,q}(X,\C)$ for $q\geq p+1$ and $(f^{n_0})^*(c_0)\leq d^{n_0}c_0$ for some strictly positive class $c_0$ in $H^{p,p}(X,\C)$. Note that the last condition implies that $\|(f^n)^*\|=O(d^n)$ on $H^{p,p}(X,\C)$.

Let $Z$ be another homogeneous compact K\"ahler manifold of dimension $m$. Fix a K\"ahler form $\omega_Z$ on $Z$ and consider on $X\times Z$ the K\"ahler form 
$\pi_X^*(\omega_X)+\pi_Z^*(\omega_Z)$, where $\pi_X,\pi_Z:X\times Z\to X, Z$ are the canonical projections.  
Let $\widehat f:X\times Z\to X\times Z$  be a bi-meromorphic map which is also a bi-holomorphic map between $X_0\times Z$ and $f(X_0)\times Z$.  Assume that $f\circ \pi_X=\pi_X\circ\widehat f$. So
$\widehat f$ preserves the vertical fibration associated with $\pi_X$. Assume finally that for $x\in X_0$, the restriction of $\widehat f$ to  $\{x\}\times Z$ is a bi-holomorphic map onto $\{f(x)\}\times Z$ whose action on Hodge cohomology is the identity. 
For the last property, we identify both $\{x\}\times Z$ and $\{f(x)\}\times Z$ to $Z$ in the canonical way. 
The property automatically holds  when $Z$ is a Grassmannian.  This is a situation we will consider later.

Here is an important proposition that we will need. Recall that Proposition \ref{prop_conj_map} allows to use this result for manifolds which are bi-meromorphic to $X\times Z$.

\begin{proposition} \label{prop_coh_prod}
Let $\delta'$ be any constant such that $\delta<\delta'<d$. Then $\|(\widehat f^n)^*\|=O(\delta'^n)$ on $H^{q,q}(X\times Z,\C)$ for $q>p+m$. There is an integer $n_1\geq 1$ such that $(\widehat f^{n_1})^*(\widehat c_1)\leq d^{n_1}\widehat c_1$ for some strictly positive class $\widehat c_1$ in $H^{p+m,p+m}(X\times Z,\C)$. 
In particular, we have $\|(\widehat f^n)^*\|=O(d^n)$ on $H^{p+m,p+m}(X\times Z,\C)$.
Moreover, if $\|(f^n)^*\|=O(\delta^n)$ on $H^{q,q}(X,\C)$ for $q<p$ then $\|(\widehat f^n)^*\|=O(d^n)$ on $H^{q,q}(X\times Z,\C)$ for every $q$. 
\end{proposition}

Observe that from the last assertion in Proposition \ref{prop_conj_map},  we only need to consider integers $n$ which are divisible by $n_0$. Therefore, replacing $f,\widehat f, d,\delta,\delta'$ with $f^{n_0},\widehat f^{n_0},d^{n_0},\delta^{n_0}$ and $\delta'^{n_0}$ allows to assume that $n_0=1$.  
We need the following lemma that can be applied to $\widehat f$ instead of $f$.

\begin{lemma} \label{lemma_pseudo-effective}
The operator $f^*$ preserves the cone of pseudo-effective classes in $H^{q,q}(X,\C)$ for every $q$. If $c$ and $c'$ are pseudo-effective classes, not necessarily of the same bi-degree,  then for every $n\geq 0$
$$(f^n)^*(c)\leq (f^*)^n(c) \qquad \mbox{and} \qquad (f^n)^*(c\smallsmile c')\leq (f^n)^*(c)\smallsmile (f^n)^*(c').$$ 
\end{lemma}
\proof
Let $c$ be a pseudo-effective class. Since $X$ is homogeneous, it can be represented by a smooth positive closed form $\alpha$. 
Observe that $f^*(\alpha)$ is a positive closed $L^1$-form which represents the class $f^*(c)$. In general, $f^*(\alpha)$ may have singularities along the indeterminacy set of $f$. So the first assertion in the lemma is clear.

Let $f':X\to X$ be another dominant meromorphic map. We first prove that $(f\circ f')^*(c)\leq f'^*(f^*(c))$ for any pseudo-effective class $c$. 
Applying inductively this inequality to $f'=f,f^2,\ldots,f^{n-1}$ gives the first inequality in the lemma.
Let $\alpha$ be a smooth positive closed  form in the class $c$. Then $(f\circ f')^*(c)$ is the class of the positive closed $L^1$-form $(f\circ f')^*(\alpha)$ and $f^*(c)$ is the class of the positive closed $L^1$-form $\beta:=f^*(\alpha)$. The above $L^1$-forms are  smooth on a suitable Zariski open subset of $X$.

Let $\beta_n$ be smooth positive forms in the class $f^*(c)$ which converge to $\beta$. Then $f'^*(\beta_n)$ are positive closed $L^1$-forms in the class $f'^*(f^*(c))$. Their masses depend only on their cohomology classes and hence are independent of $n$. Extracting a subsequence allows to assume that $f'^*(\beta_n)$ converge to some positive closed current $\gamma$. We can obtain $\beta_n$ from $\beta$ using a convolution with holomorphic automorphisms close to the identity as mentioned above. We get that $\beta_n$ converge to $\beta$ locally uniformly outside the singularities of $\beta$. It follows that $\gamma$ is equal to $f'^*(\beta)$ on a Zariski open subset of $X$. Since $f'^*(\beta)$ is equal to the $L^1$-form $(f\circ f')^*(\alpha)$ on a Zariski open set, we deduce that $\gamma\geq (f\circ f')^*(\alpha)$ since $\gamma$ may have a singular part supported by a subvariety of $X$. Thus $(f\circ f')^*(c)\leq f'^*(f^*(c))$. 
This implies the first inequality in the lemma.

For the last inequality in the lemma, we can for simplicity assume that $n=1$. Consider  a smooth positive closed form $\theta$ in the class $c'$. So $f^*(\alpha\wedge\theta)$ represents the class $f^*(c\smallsmile c')$ and $\beta_n\wedge  f^*(\theta)$ represents $f^*(c)\smallsmile f^*(c')$. So any limit value of the sequence $\beta_n\wedge f^*(\theta)$  represents $f^*(c)\smallsmile  f^*(c')$. Such a limit value
is equal to $ f^*(\alpha)\wedge f^*(\theta)$ on a Zariski open set and hence equal to $f^*(\alpha\wedge\theta)$ on a Zariski open set. The last current is an $L^1$-form and has no mass on proper analytic subsets of $X$. We conclude that the considered  limit value is at least equal to $f^*(\alpha\wedge \theta)$. The last inequality in the lemma follows.
\endproof 
 
Recall that we assumed $n_0=1$. So the last lemma implies that $(f^n)^*(c_0)\leq d^nc_0$ for every $n\geq 1$. 
Denote by $c_X$ and $c_Z$ the classes of $\omega_X$ and $\omega_Z$. Their powers  are strictly positive classes. Fix a constant $\delta_1$ such that $\delta<\delta_1<\delta'$.
Replacing $f$ with a power of $f$ allows to assume that $f^*(c_X^q)\leq \delta_1c_X^q$ for $q>p$. 
We also have the same inequality for $q<p$ when we assume that $\|(f^n)^*\|=O(\delta^n)$ on $H^{q,q}(X,\C)$ for $q<p$.

Recall that we assume that $H^{r,s}(X,\C)=0$ for $r\not=s$.  By K\"unneth formula, we have the following canonical decomposition of the Hodge cohomology on $X\times Z$
$$H^{q,q}(X\times Z,\C)=\sum_r H^{r,r}(X,\C)\times H^{q-r,q-r}(Z,\C).$$
In general, the above decomposition is not invariant under the action of $\widehat f^*$. 
Define 
$$E_{q,s}:=\sum_{r\geq s} H^{r,r}(X,\C)\times H^{q-r,q-r}(Z,\C).$$
So we have a decreasing sequence of vector spaces with $E_{q,0}=H^{q,q}(X\times Z,\C)$.
We can show that $E_{q,s}$ is invariant under $\widehat f^*$. 
 
 \begin{lemma} \label{lemma_filtre}
 There is a constant $A>0$ such that 
 $$\widehat f^*\big(c_X^s\otimes c_Z^{q-s}\big)\leq \delta_1 (c_X^s\otimes c_Z^{q-s}) + A \sum_{j=1}^{q-s} c_X^{s+j}\otimes c_Z^{q-s-j}$$
 for all $p+1\leq s\leq k$ and $0\leq q-s\leq m$.
  If $f^*(c_X^s)\leq\delta_1 c_X^s$ for $s<p$, then the above estimate also holds for $s< p$.
  \end{lemma}
 \proof
We prove the first assertion. The second one is obtained in the same way. 
We first show that $\widehat f^*\big(c_X^s\otimes c_Z^{q-s}\big)$ belongs to $E_{q,s}$. 
For this purpose, it is enough to check that the cup-product of this class with $\pi_X^*(c)$ vanishes for any class $c$ in $H^{k-s+1,k-s+1}(X,\C)$. Moreover, we only have to check the last property for $c$ pseudo-effective because such classes generate $H^{k-s+1,k-s+1}(X,\C)$. 
Since the considered class is pseudo-effective, we only need to show that its cup-product with $\pi_X^*(c)$ is negative or zero.

 Define $\widetilde c_X:=\pi_X^*(c_X)$ and $\widetilde c_Z:=\pi_Z^*(c_Z)$. By Lemma \ref{lemma_pseudo-effective}, we have
 \begin{eqnarray*}
 \widehat f^*\big(c_X^s\otimes c_Z^{q-s}\big) &\leq &  \widehat f^*(\widetilde c_X^s)\smallsmile \widehat f^*(\widetilde c_Z^{q-s})=\pi_X^*(f^*(c_X^s))\smallsmile \widehat f^*(\widetilde c_Z^{q-s})\\
 &\leq & \delta_1 \pi_X^*(c_X^s)\smallsmile \widehat f^*(c_Z^{q-s}).
 \end{eqnarray*}
The cup-product of the factor $\pi_X^*(c_X^s)$ with $\pi_X^*(c)$ vanishes for a bi-degree reason.
It follows that   $\widehat f^*\big(c_X^s\otimes c_Z^{q-s}\big)$ belongs to $E_{q,s}$.   
  
Since $c_X$ and $c_Z$ are K\"ahler classes, any real class in $E_{q,s+1}$ can be bounded by the last term of the inequality in the lemma provided that $A$ is large enough. Therefore, 
in order to obtain the result, we only need to check that  $\pi_X^*(c_X^s)\smallsmile \widehat f^*(c_Z^{q-s})$ is equal to $c_X^s\otimes c_Z^{q-s}$ plus a class in $E_{q,s+1}$. By Poincar\'e's duality, if $\kappa_X$ and $\kappa_Z$ are classes in $H^{k-s,k-s}(X,\C)$ and $H^{m-q+s,m-q+s}(Z,\C)$ respectively with $c_X^s\smallsmile \kappa_X=1$ and $c_Z^{q-s}\smallsmile \kappa_Z=1$, it suffices to show that 
$$\pi_X^*(c_X^s)\smallsmile \widehat f^*(c_Z^{q-s})\smallsmile \pi_X^*(\kappa_X)\smallsmile \pi_Z^*(\kappa_Z)=1.$$
But this identity is clear because $\pi_X^*(c_X^s)\smallsmile\pi_X^*(\kappa_X)$ can be represented by a generic fiber of $\pi_X$ and the restriction of $\widehat f$ to generic fibers of $\pi_X$ acts trivially on the cohomology of these fibers. The lemma follows.  
\endproof
 
 Using exactly the same arguments gives us the following lemma.

 \begin{lemma} \label{lemma_filtre_bis}
 There is a constant $A>0$ such that 
 $$\widehat f^*\big(c_X^p\otimes c_Z^{q-p}\big)\leq d (c_X^p\otimes c_Z^{q-p}) + A \sum_{j=1}^{q-p} c_X^{p+j}\otimes c_Z^{q-p-j}$$
 for  $0\leq q-p\leq m$. 
  \end{lemma}
 
\noindent
{\bf End of the proof of Proposition \ref{prop_coh_prod}.} 
Recall that we assumed  $n_0=0$ and $f^*(c_X^s)\leq \delta_1 c_X^s$ for $s>p$. For the last assertion in the proposition, we can also assume that $f^*(c_X^s)\leq \delta_1 c_X^s$ for $s<p$.

\medskip
\noindent
{\bf Claim.}
We have $\|(\widehat f^n)^*(c_X^s\otimes c_Z^{q-s})\|=O(d^n)$ for $s=p$ and $\|(\widehat f^n)^*(c_X^s\otimes c_Z^{q-s})\|=O(n^M\delta_1^n)$ for $p+1\leq s\leq k$ and $0\leq q-s\leq m$ with a suitable integer $M$.  For the last assertion in the proposition, we have  $\|(\widehat f^n)^*(c_X^s\otimes c_Z^{q-s})\|=O(d^n)$ for every $s$.

\medskip

It is not difficult to deduce the first and third assertions in the proposition from the claim and K\"unneth decomposition of cohomology on $X\times Z$. 
For the first assertion, we use that 
$H^{q,q}(X\times Z,\C)=E_{q,p+1}$ when $q>p+m$. 

We prove now the claim. By Lemma \ref{lemma_pseudo-effective}, it is enough to check the same estimates for $(\widehat f^*)^n$ instead of $(\widehat f^n)^*$. 
We will obtain these estimates using a decreasing induction on $s$. 
We can consider that the case $s=k+1$ is trivial because $c_X^{k+1}=0$ for a bi-degree reason. So assume that our above claim is true for $s+1,s+2,\ldots$ instead of $s$ with $0\leq s\leq k$ and we prove it for $s$. We only consider the case $s\leq p-1$ and the other cases can be obtained in the same way. By Lemma \ref{lemma_filtre}, we have
$$(\widehat f^*)^N(c_X^s\otimes c_Z^{q-s})\leq \delta_1 (\widehat f^*)^{N-1}\big(c_X^s\otimes c_Z^{q-s}\big) + 
A \sum_{j=1}^{q-s} (\widehat f^*)^{N-1}\big(c_X^{s+j}\otimes c_Z^{q-s-j}\big).$$
By induction hypothesis, the norm of the last sum is smaller than a constant times $d^{N-1}$. Since $\delta_1<d$, the last inequality applied to $N=n,n-1,\ldots,1$ implies that $\|(\widehat f^*)^n(c_X^s\otimes c_Z^{q-s})\|=O(d^n)$. The claim follows.

It remains to prove the second assertion in the proposition. Observe first that $H^{p+m,p+m}(X\times Z,\C)=E_{p+m,p}$. Therefore, any combination of $c_X^s\otimes c_Z^{p+m-s}$ with strictly positive coefficients and with $p\leq s\leq \min(k,p+m)$ is a strictly positive class in  $H^{p+m,p+m}(X\times Z,\C)$. Define
$$\widehat c_1:=c_X^p\otimes c_Z^m+\sum_{j=1}^{\min(k-p,m)} A_jc_X^{p+j}\otimes c_Z^{m-j},$$
where $A_j$ are constants large enough such that $A_j\ll A_{j+1}$. It is not difficult to deduce from Lemmas \ref{lemma_filtre} and \ref{lemma_filtre_bis} that $\widehat c_1$ satisfies the proposition. Here, we can take $n_1=1$ but we already replaced twice $f$ with an iterate. 
\hfill $\square$

\bigskip
\noindent
{\bf Main examples.} 
We describe now the main examples that will be considered later. 
Let $X$ be a compact K\"ahler manifold of dimension $k$ which is a projective space or the product of two projective spaces. So $X$ is homogeneous and $H^{r,s}(X,\C)=0$ for $r\not=s$. Let $\omega_X$ be a K\"ahler form on $X$ and denote by $c_X$ its cohomology class. 
Let $f:X\to X$ be a bi-rational map and let $d,p,\delta$ be as above such that $f^*(c_X^p)\leq dc_X^p$ and $f^*(c_X^q)\leq\delta c_X^q$ for $q\not=p$.

Denote by $\Gr(X,k-p)$ the  space of points $(x,[v])$, where $x$ is a point in $X$ and $[v]$ is the direction of a {\it simple} complex tangent $(k-p)$-vector of $X$ at $x$.
The natural projection from $\Gr(X,k-p)$ to $X$ defines a fibration whose fibers are isomorphic to the Grassmannian $\G$ of vector subspaces of dimension $k-p$ in $\C^k$. We can lift $f$ to a bi-rational map $\widehat f:\Gr(X,k-p)\to\Gr(X,k-p)$ by $\widehat f(x,[v]):=(f(x), [f_*(v)])$ for $x$ in a suitable Zariski open subset of $X$.  Let $\widetilde f$ denote the lift of $\widehat f$ to $\Gr(\Gr(X,k-p),k-p)$ which is defined in the same way.

\begin{proposition} \label{prop_ex}
We have $\|(\widehat f^n)^*\|=O(d^n)$  on the Hodge cohomology groups of $\Gr(X,k-p)$. If $K$ is the dimension of 
$\Gr(\Gr(\P^k,k-p),k-p)$, we also have  $\|(\widetilde f^n)^*\|=O(d^n)$ on the Hodge cohomology group 
of $\Gr(\Gr(\P^k,k-p),k-p)$ of bi-degree $(q,q)$ for every $q\geq K-k+p$. 
\end{proposition}
\proof
Over a chart $\C^k\subset X$, the fibration $\Gr(X,k-p)$ can be identified in a natural way with $\C^k\times\G$. So $\Gr(X,k-p)$ is bi-rational to $X\times\G$. By Proposition \ref{prop_conj_map},  we can consider $\widehat f$ as a bi-rational map of $X\times\G$. The first assertion in the proposition is a direct consequence of the last assertion in Proposition \ref{prop_coh_prod}. 
The second assertion in Proposition \ref{prop_coh_prod} can be applied to $f':=\widehat f$ and $X':=X\times \G$ as we will see below.

The manifold $\Gr(\Gr(\P^k,k-p),k-p)$ is bi-rational to $X'\times\G'$, where $\G'$ is a Grassmannian. In order to see this point, it is enough to identify some Zariski open subset of $X'$ with a Zariski open subset of a complex Euclidean space. So we can consider 
$\widetilde f$ as a map on $X'\times \G'$ which preserves the natural fibration over $X'$. 
Applying Proposition \ref{prop_coh_prod} to $f',\widetilde f, X'\times \G'$ instead of $f,\widehat f, X\times Z$ gives the result.
\endproof

\begin{examples} \rm \label{ex_Henon}
Let $f$ be a H\'enon automorphism on $\C^k$ that we extend to a bi-rational map on $\P^k$. Let $d_\pm, d, p$ be as in the introduction. 
So the operator $f^*$ on $H^{q,q}(\P^k,\C)$ is just the multiplication by $d_+^q$ for $q\leq p$ and by $d_-^{k-q}$ for $q\geq p$. 
Define $\delta:=\max(d_+^{p-1},d_-^{k-p-1})$. We have $1\leq \delta< d$. So the last proposition can be applied to $f$ and its lifts to $\Gr(\P^k,k-p)$ and $\Gr(\Gr(\P^k,k-p),k-p)$. 

Consider the bi-rational map $F=(f,f^{-1})$ on $\P^k\times \P^k$. 
If $\omega_\FS$ denotes the Fubini-Study form on $\P^k$ and $\pi_1,\pi_2$ denote the projections from $\P^k\times\P^k$ onto its factors, we consider on $\P^k\times\P^k$ the K\"ahler metric $\pi_1^*(\omega_\FS)+\pi_2^*(\omega_\FS)$. 
Let $c_\FS$ denote the class of 
$\omega_\FS$. We have $F^*(c_\FS^p\otimes c_\FS^{k-p})=d^2 (c_\FS^p\otimes c_\FS^{k-p})$ and 
$F^*(c_\FS^r\otimes c_\FS^s)\leq d\delta (c_\FS^r\otimes c_\FS^s)$ for $(r,s)\not=(p,k-p)$. So we can apply Proposition \ref{prop_ex} to $F,d^2,d\delta,k$ instead of $f,d,\delta,p$
\end{examples}

\section{Woven currents and tame currents} \label{section_woven}

This section contains some geometric properties of positive closed currents that we will  use in our study of H\'enon maps. We discuss the notions of woven and tame currents which have an independent interest. 
Laminar currents in dimension 2 were introduced and studied by  Bedford-Lyubich-Smillie \cite{BLS}, see also Sullivan \cite{Sullivan}. Woven currents and laminar currents in higher dimension were introduced by the first author of the present paper in \cite{Dinh}. 

If two Riemann surfaces in $\C^2$ are the limits of two sequences of Riemann surfaces $\Gamma_n$ and $\Gamma_n'$ with $\Gamma_n\cap\Gamma_n'=\varnothing$, then their intersection is either empty or also a  Riemann surface. In higher dimension and codimension, this property is no longer true. This is one of the main difficulties  with woven and laminar currents in higher dimension. In order to simplify the exposition, we will work on a projective manifold $V$ of dimension $l$. We can also extend the theory to currents on non-compact manifolds. Fix a K\"ahler form $\omega$ on $V$ and consider the K\"ahler metric on $V$ induced by $\omega$.

\medskip\noindent
{\bf Measurable webs, woven currents and standard refinement.}
We introduce here some basic notions. 
Let $0\leq r\leq l$ be an integer. Denote by $\Lam_r(V)$ the set of 
positive $(l-r,l-r)$-currents that can be written as a finite or countable sum $S=\sum [\Lambda_i]$, where the $\Lambda_i$'s are irreducible analytic sets of dimension $r$, immersed in $V$,  such that $\sum\volume(\Lambda_i)<\infty$. Here $[\Lambda_i]$ denotes the current of integration on $\Lambda_i$ and $\volume(\Lambda_i)$ denotes the $2r$-dimensional volume of $\Lambda_i$ that can be computed using  Wirtinger's theorem by 
$$\volume (\Lambda_i)={1\over r!} \int_{\Lambda_i} \omega^r.$$
We say that $S$ is a {\it lame}.
Note that we don't assume that the $\Lambda_i$'s are disjoint.
Moreover, given a current $S$ as above the decomposition  $S=\sum [\Lambda_i]$ is not unique. The reason to consider here finite or countable sums is to gain flexibility in working with woven and laminar currents.

We can identify $\Lam_r(V)$ with a subset of the family $\Pc_r(V)$ of (strongly) positive currents of bi-dimension $(r,r)$ on $V$. The later is 
a metric space endowed with the following family of distances
$$\dist_\alpha(T,T'):=\sup_{\|\varphi\|_{\Cc^\alpha}\leq 1} |\langle T-T',\varphi\rangle|.$$
The induced topology on $\Pc_r(V)$ is the same for any $\alpha>0$ and coincides with the weak topology on currents. 
Define $\Lam_r^*(V):=\Lam_r(V)\setminus\{0\}$ and 
$\Pc_r^*(V):=\Pc_r(V)\setminus\{0\}$. Using the local description below of woven currents, we will see in Lemma \ref{lemma_Borel} that $\Lam_r^*(V)$ is a universally measurable set, i.e. it is measurable with respect to all Borel probability measures on $\Pc_r(V)$. 

\begin{definition} \rm
We call {\it measurable $l$-web} any positive Borel measure $\nu$ on $\Pc_r^*(V)$ without mass outside $\Lam_r^*(V)$ such that 
$$\int_{\Lam_r^*(V)} \|S\| d\nu(S)<\infty.$$
A  current $T$ of bi-dimension $(r,r)$ on $V$ is said to be {\it woven} if there is a measurable $r$-web $\nu$ (which is called a measurable $r$-web  {\it associated} with $T$) such that 
$$T:= \int_{\Lam_r^*(V)} S d\nu(S)$$
or equivalently 
$$\langle T,\varphi\rangle := \int_{\Lam_r^*(V)} \langle S,\varphi\rangle  d\nu(S)$$
for any test continuous $(r,r)$-form $\varphi$ on $V$. 
\end{definition}

Note that a woven current may be associated with different measurable webs, e.g. the Fubini-Study form $\omega_\FS$ on $\P^k$ can be obtained as an average of hypersurfaces of degree $d$ for any positive integer $d$ as shown in the following example.

\begin{example} \label{ex_Fubini}\rm
Let $\U(k+1)$ denote the unitary group which acts naturally on $\P^k$. If $\sigma$ is the Haar measure on $\U(k+1)$ and $H$ is a subvariety of dimension $r$ and of degree $d$ of $\P^k$, we have the following identity in the sense of currents
$$\omega_\FS^{k-r}=d^{-1}\int_{\tau\in \U(k+1)} \tau_*[H]d\sigma(\tau).$$
This identity says that $\omega_\FS^{k-r}$ can be written as an average of currents of integration on  subvarieties of degree $d$ and of dimension $r$. So it is a woven current. 
\end{example}

Note  that lames can be  divided into smaller ones giving  different webs associated with the same current.
If $\nu$ is a measurable $r$-web, it has finite mass outside any neighbourhood of 0 in $\Pc_r(V)$. Indeed, outside any neighbourhood of 0, currents in $\Pc_r(V)$ have mass bounded from below by a strictly positive constant.

Let $\Sigma_r$ be the set of $(R,S)$ in $\Lam^*_r(V)\times\Lam^*_r(V)$ such that $R\leq S$. Denote by $\pi$ and $\pi'$ the natural projections $(R,S)\mapsto R$ and $(R,S)\mapsto S$. 

\begin{definition} \rm \label{def_web_refinement}
Let $\nu_1$ and $\nu_2$ be two measurable $r$-webs on $V$. 
We say that $\nu_1$ is {\it a refinement} of $\nu_2$ and  write $\nu_2\prec \nu_1$
if there is a positive measure $\nu_{12}$ on $\Sigma_r$  such that
\begin{enumerate}
\item $\nu_1=\pi_*(\nu_{12})$. 
\item If $\nu_{12}'$ is the restriction of $\nu_{12}$ to the complement of a neighbourhood of 0 in $\Lam_r(V)$, then $\pi'_*(\nu_{12}')$ is absolutely continuous with respect to $\nu_2$.
\item For $\nu_2$-almost every $S$
denote by
$\nu_{12}^S:=\langle \nu_{12}|\pi'|S\rangle$  the conditional measure of $\nu_{12}$ with respect to the fiber of $\pi'$ over the point $S$. We identify it to a measure on $\Lam_r^*(V)$. Then $\nu_{12}^S$  is a finite or countable sum of Dirac masses and defines a measurable $l$-web associated with $S$.
\end{enumerate}
Two measurable $r$-webs are {\it equivalent} if they admit a common refinement, see also Lemma \ref{lemma_web_eq} below.
\end{definition}

Roughly speaking, in order to get a refinement $\nu_1$ of   $\nu_2$ we decompose the lames $S$ of $\nu_2$ into a finite or countable number of smaller ones using the conditional measures $\nu_{12}^S$. For example, when $\nu_2$ is a Dirac mass at a point $R$, then $\nu_{12}^S=0$ for $S\not=R$ and $\nu_{12}^R=\nu_{12}$. In this case, $\nu_{12}^R$ is identified with $\nu_1$. We get a decomposition of $R$ into a finite or countable sum of lames. The general case can be deduced from this case by taking an average with respect to $\nu_2$. 

\begin{lemma}
Let $\nu_1$ and $\nu_2$ be measurable $r$-webs on $V$. Assume that $\nu_2\prec\nu_1$. Then they are associated with the same current.
If $\nu_3$ is another measurable $r$-web such that  $\nu_3\prec\nu_2$, then $\nu_3\prec\nu_1$. 
\end{lemma}
\proof
The properties in Definition \ref{def_web_refinement} imply that 
$$\int_R R d\nu_1(R)=\int_{R,S} R d\nu_{12}(R,S)=\int_S \Big(\int_R R d\nu_{12}^S(R)\Big) d\nu_2(S)=\int_S Sd\nu_2(S).$$
This gives the first assertion in the lemma.

Denote by $\nu_{12}$ the measure  in Definition \ref{def_web_refinement} and $\nu_{23}$ the similar one associated with $\nu_2$ and $\nu_3$. If we identify each fiber $\{S\}\times\Lam_r^*(V)$ of $\pi'$ with $\Lam_r^*(V)$, we can define a measure $\nu_{13}$ on $\Lam_r^*(V)\times\Lam_r^*(V)$ by their conditional measures with respect to $\nu_3$
$$\langle\nu_{13}|\pi'|S\rangle:=\int_R \nu_{12}^R d \nu_{23}^S(R).$$
It is not difficult to check that $\nu_{13}$ induces the relation $\nu_3\prec\nu_1$. Indeed, one can consider the case where $\nu_3$ is a Dirac mass and obtain the general case by taking an average. 
\endproof

The abstract formalism on woven currents introduced above does not require a choice of local coordinates. It does not depend on the metric on the manifold. So it offers  a convenient setting to work with different local coordinate systems and other operations on woven currents, e.g. the lifting of currents to Grassmannian bundles. However,
in order to get a more precise picture on woven currents and to construct measurable webs, 
we will work in convenient local coordinates. 
Up to a choice of local coordinate systems, we give now a uniform way to decompose a measurable web into an infinite sum of $s$-elementary webs  with $s=1,2,\ldots$

\medskip\noindent
\underline{\it Approximation and refinement for local analytic sets.} Let $z=(z_1,\ldots,z_l)$ be local holomorphic coordinates with $|z_i|<4$. Write $z_i=x_i+\sqrt{-1} y_i$ with $x_i,y_i\in\R$. We consider here measurable webs $\nu$ with lames inside the cube $\U:=\{|x_i|,|y_i|\leq 1\}=[-1,1]^{2l}$.   
Let $\Lambda$ be an irreducible analytic set of dimension $r$, not necessarily closed,  immersed in the cube and with finite $2r$-dimensional volume. Then there is a projection onto $r$ coordinates among $z_1,\ldots, z_l$ whose generic fibers intersect  $\Lambda$ in  finite or countable sets. For simplicity, assume that
 the projection $z\mapsto (z_1,\ldots,z_r)$ satisfies this property. Otherwise, in order to get the same construction for all varieties, we choose the convenient projection with the smallest lexicographical index. 

For each $s\geq 1$, consider the real hyperplanes $\{x_i=j2^{-s}\}$ and $\{y_i=j2^{-s}\}$ with $j\in\Z$ and $i\leq r$ that will be called {\it separating hyperplanes} (there are also separating hyperplanes associated with the other projections). They divide the cube $[-2,2]^{2l}$ into $2^{2r(s+2)}$ (closed) towers. 
Two towers are said to be {\it adjacent} if they have common points. In particular, a tower is adjacent to itself. The union of the towers which are adjacent to a given tower is called a {\it fat tower}. The union of towers which are adjacent to one of the towers in the last fat tower is called a {\it very fat tower}. So fat and very fat towers are just extensions in the horizontal directions of ordinary towers. 

Consider complex manifolds  which are graphs over the basis of a tower that can be extended to graphs over the basis of the associated very fat tower (we don't ask the later graphs to be defined over the boundary of the basis of the very fat tower). We call them {\it $s$-elementary lames}. 
The extension of a graph in a tower to a graph in the corresponding fat or very fat tower is called the {\it fat or very fat extension} respectively. We will work with the family of $s$-elementary lames in the box $\U:=[-1,1]^{2l}$. 
Observe that this family of  is compact with respect to the local uniform convergence topology on holomorphic graphs.  
A finite positive measure on this compact set is called an {\it $s$-elementary web} in $\U$. 
Considering the fat and very fat  extensions of graphs is just a technical point which allows us to avoid the possible bad behavior of graphs near the vertical boundary of a tower.

Denote by $\nu^{[s]}_\Lambda$ the maximal  $s$-elementary web whose associated current is smaller than or equal to  $[\Lambda]$. 
Denote also by $[\Lambda^{[s]}]$ this current. 
The measure 
$\nu^{[s]}_\Lambda$ is a sum of Dirac masses. It  gives  an approximation of $\Lambda$: we have $[\Lambda]-[\Lambda^{[s]}]\to 0$ in the mass norm as $s\to\infty$. Note that $[\Lambda]-[\Lambda^{[s]}]$ is a lame, i.e. a point in $\Lam_r(V)$.
If it is not zero, the sum of the Dirac mass at this point and $\nu_\Lambda^{[s]}$ is denoted by $\nu_\Lambda^{(s)}$. Otherwise, define $\nu_\Lambda^{(s)}:=\nu_\Lambda^{[s]}$. So $\nu_\Lambda^{(s)}$ is a measure associated with $[\Lambda]$ and gives a refinement of $[\Lambda]$. Finally, we observe that the construction depends only on the current $[\Lambda]$. More precisely, if we remove from $\Lambda$ a closed subset of zero $2r$-dimensional measure, we obtain another analytic set  associated with the same current. Our construction gives the same approximation and refinement.

\medskip\noindent
\underline{\it Approximation and refinement for local lames.}
The construction extends without difficulty to any lame $S=\sum [\Lambda_i]$ in $\U$
 such that the above condition on the projection 
$z\mapsto(z_1,\ldots, z_r)$ is satisfied for each component $\Lambda_i$. 
We define $\nu^{[s]}_S$ as the maximal $s$-elementary web whose associated current, denoted by $S^{[s]}$, is smaller than or equal to $S$. The current $S-S^{[s]}$ is a lame, i.e. a point in $\Lam_r(V)$. If it is not zero, the sum of the Dirac mass at this point and $\nu_S^{[s]}$ is denoted by $\nu_S^{(s)}$. Otherwise, define $\nu_S^{(s)}:=\nu_S^{[s]}$. This is a refinement of $S$.

We extend the construction to an arbitrary lame $S=\sum [\Lambda_i]$ in $\U$.
Denote by $S_1$ the sum of $[\Lambda_i]$ which satisfy the above condition on the projection $z\mapsto(z_1,\ldots, z_r)$. The above construction can be applied to $S_1$. Denote by $S_2$ the sum of the remaining $[\Lambda_i]$ satisfying the similar condition for the next projection $z\mapsto (z_{i_1},\ldots,z_{i_r})$ with respect to the lexicographical index order. We do the similar construction for the new projection and repeat it again for the other projections respecting always the lexicographical order. With the notations similar to the ones given above, define
$$\nu_S^{[s]}:=\sum_i \nu_{S_i}^{[s]} \quad \text{and} \quad \nu_S^{(s)}:=\sum_i \nu_{S_i}^{(s)}.$$
They are respectively the {\it $s$-approximation} and {\it $s$-refinement} of $S$. They do not depend on the choice of the decomposition $S=\sum [\Lambda_i]$.

\medskip\noindent
\underline{\it Approximation and refinement for global lames.}
Fix a covering of  $V$ by a finite number of cubes $\U_1,\ldots,\U_N$ as above.
Assume that $S$ is a lame in $V$. We can decompose it into local lames 
$S=S_1+\cdots+S_N$, where $S_1$ is the restriction of $S$ to $\U_1$ and by induction $S_i$ is the restriction to $\U_i$ of $S-S_1-\cdots-S_{i-1}$. We then apply the above construction to each $S_i$ in $\U_i$ and obtain the webs $\nu^{(s)}_{S_i}$ and $\nu^{[s]}_{S_i}$. Define 
$$\nu^{(s)}_S:=\nu^{(s)}_{S_1}+\cdots +\nu^{(s)}_{S_N} \quad \text{and}\quad \nu^{[s]}_S:=\nu^{[s]}_{S_1}+\cdots +\nu^{[s]}_{S_N}.$$ 
The sum of $N$ $s$-elementary webs  on $\U_1,\ldots,\U_N$ respectively is called  an {\it $s$-elementary web} on $V$. If $S^{[s]}$ is the current associated with $\nu_S^{[s]}$, then $S-S^{[s]}\to 0$ in the mass norm.

\medskip\noindent
\underline{\it Approximation and  refinement for global webs.}
We can apply the same method to refine or  approximate all woven currents. If $T$ is such a current and $\nu$ is an associated web, define using the above notations 
$$\nu^{(s)}:=\int_{\Lam_r^*(V)} \nu^{(s)}_S d\nu(S) \quad \text{and}\quad \nu^{[s]}:=\int_{\Lam_r^*(V)} \nu^{[s]}_S d\nu(S).$$
We say that  $\nu^{(s)}$ and $\nu^{[s]}$ are respectively the {\it standard $s$-refinement} and the {\it standard $s$-approximation} of $\nu$.  If $T^{[s]}$ is the current associated with $\nu_T^{[s]}$, then $T-T^{[s]}\to 0$ in the mass norm. 
Here are two applications of the construction.

\begin{lemma} \label{lemma_Borel}
The set  $\Lam_r^*(V)$ is universally measurable. 
\end{lemma}
\proof
Using a finite number of boxes $\U_1,\ldots,\U_N$ as above, it is not difficult to reduce the problem to the set of  
lames $S$ with 
supports in a box $\overline \U$ such that the projection $\pi(z):=(z_1,\ldots, z_r)$ is of maximal rank on each component of $S$. Denote by $\Lc$ the set of such lames $S$ with mass  bounded by a fixed constant $M$. It is enough to show that $\Lc$ is universally measurable. 

Using the decomposition of lames into elementary ones as above, we see that  $\Lc$ is also the set of currents of the forms $S=\sum_{s\geq 1} S_s$, where $S_s$ is a finite sum of $s$-elementary lames and $\|S\|\leq M$. Note that such a current $S$ is associated with infinitely many 
different decompositions into elementary lames.  

Denote by $\Lc_s$ the set of currents which are equal to a sum of at most $2^{2rs}M$ $s$-elementary lames. This is a compact set of currents.  Consider the infinite product space $\Pi_{s\geq 1} \Lc_s$ endowed with the natural product topology. Let $\Lc'$ be the subset of points $(S_1,S_2,\ldots)$ in this space such that $\|S\|\leq M$. 
The last condition is equivalent to $\sum_{s=1}^n\|S_s\|\leq M$ for every $n\geq 1$. 
So $\Lc'$ is a Borel set and the 
map $(S_1,S_2,\ldots)\mapsto S:=\sum S_s$ from $\Lc'$ to $\Lc$ is continuous and surjective. Therefore, the image $\Lc$ of this map is universally measurable, see \cite[p.98]{DM}. The lemma follows.
\endproof

\begin{lemma} \label{lemma_web_eq}
Let $\nu_1,\nu_2$ and $\nu_3$ be measurable $r$-webs on $V$.
If $\nu_3\prec \nu_1$ and $\nu_3\prec\nu_2$, then there is a measurable $r$-web $\nu_0$ such that $\nu_1\prec\nu_0$ and $\nu_2\prec\nu_0$. 
In particular, if 
$\nu_1$ and $\nu_2$ are equivalent to $\nu_3$ in the sense of Definition \ref{def_web_refinement},  then $\nu_1$ and $\nu_2$ are also equivalent. 
\end{lemma}
\proof
Using a covering of $V$ by cubes $\U_1,\ldots, \U_N$ as above, we can reduce the problem to the case where 
 $T$ is a current on the cube $[-1,1]^{2l}$ of $\C^l$.
Using the above construction with the similar notations,
we have 
$\nu_3^{(s)}\prec \nu_1^{(s)}$, $\nu_3^{(s)}\prec \nu_2^{(s)}$,
$\nu_1^{[s]}\leq \nu_3^{[s]}$ and $\nu_2^{[s]}\leq \nu_3^{[s]}$. The last two inequalities are not in general equalities because 
some $s$-elementary lames in $\nu_3^{[s]}$ may not be  lames of $\nu_1^{[s]}$ and $\nu_2^{[s]}$ after the refinement.
The measurable webs $\nu_i^{(s)}$ are all associated with $T$. 
If $T_i^{[s]}$ is the woven current associated with $\nu_i^{[s]}$, then $T-T_i^{[s]}$ is a woven current which tends to 0 as $s\to\infty$.
Write $ \nu_i^{[s]}=h_i \nu_3^{[s]}$ with $0\leq h_i\leq 1$. Define $h:=\min(h_1,h_2)$, $\vartheta^{[s]}:=h\nu_3^{[s]}$ and denote by $T^{[s]}$ the woven current associated with $\vartheta^{[s]}$. Since $1-h\leq (1-h_1)+(1-h_2)$, the woven current $T-T^{[s]}$  tends to 0 as $s\to\infty$. 

For a fixed integer $s_0$ large enough, define $\vartheta\langle 1\rangle :=\vartheta^{[s_0]}$ and $\nu_i\langle 1\rangle :=\nu^{(s_0)}_i-\vartheta^{[s_0]}$ for $i=1,2,3$. 
We can see $\nu_i^{(s_0)}$ as refinements of $\nu_i$ satisfying the same hypotheses on $\nu_i$. The web $\vartheta\langle 1\rangle $ is approximately a common refinement of $\nu_i$.
With $s_0$ large enough, the mass of the woven current associated with $\nu_i\langle 1\rangle $, i.e. the error of the approximation, which does not depend on $i$, is smaller than $1/2$.
We repeat the above construction in order to refine approximately $\nu_i\langle 1\rangle $. We obtain a web $\vartheta\langle 2\rangle $ such that the woven current associated with $\nu_i\langle 1\rangle $ is approched by  the one associated with $\vartheta\langle 2\rangle $: the difference of these currents is associated with three measurable webs $\nu_i\langle 2\rangle$ and has a mass
smaller than $1/4$. 

By induction, we obtain sequences $\nu_i\langle m\rangle$ and $\vartheta\langle m\rangle$ such that 
the mass of the woven current associated with $\nu_i\langle m\rangle$ is smaller than $2^{-m}$. By construction, the web $\nu_0:=\sum \vartheta\langle m\rangle$ refines all $\nu_1$, $\nu_2$ and $\nu_3$. This completes the proof of the lemma.
\endproof

\noindent
{\bf Weakly laminar, laminar and tame currents.} We introduce now currents with stronger geometric properties. 
\begin{definition} \rm
A woven current $T$ of bi-dimension $(r,r)$ on $V$ is {\it weakly laminar}\footnote{The terminology is changed with respect to the one in \cite{Dinh}.} if it admits a measurable $r$-web $\nu$ such that for $\nu\times\nu$-almost every pair of lames $S=\sum [\Lambda_i]$ and $S'=\sum [\Lambda_j']$ either $\Lambda_i\cap \Lambda_j'=\varnothing$ or $\Lambda_i\cap \Lambda_j'$ is open in $\Lambda_i$ and in $\Lambda_j'$ for all $i,j$. We say that $T$ is {\it laminar} if there is a measurable $r$-web $\nu$ associated with $T$ and a measurable subset $A$ of $\Lam_r^*(V)$ such that $\nu=0$ outside $A$ and for all $S=\sum [\Lambda_i]$ and $S'=\sum [\Lambda_j']$ in $A$ either $\Lambda_i\cap \Lambda_j'=\varnothing$ or $\Lambda_i\cap \Lambda_j'$ is open in $\Lambda_i$ and in $\Lambda_j'$ for all $i,j$. We say that $T$ is {\it completely weakly laminar} or {\it completely laminar} if the above corresponding property holds for all measurable webs $\nu$ associated with $T$.  
\end{definition}

A priori, the pairs of lames satisfying the condition for the weak laminarity property form a subset in $\Lam_r^*(V)\times \Lam_r^*(V)$ which is not necessarily of the product form $A\times A$. So laminar currents are weakly laminar. 
One can show that the converse holds when $r=\dim V-1$. 
If $\Lambda_1$ and $\Lambda_2$ are two manifolds of dimension $r$ such that $\Lambda_1\cap \Lambda_2$ is non-empty and of dimension $<r$ then $[\Lambda_1]+[\Lambda_2]$ is laminar but it is not completely weakly laminar. To see this point, we can consider the web which is the sum of the Dirac masses at $[\Lambda_1]$ and at $[\Lambda_2\setminus\Lambda_1]$. 
We have the following proposition which was obtained by Dujardin for $(1,1)$-currents on manifolds of dimension 2 \cite{Dujardin2}. 

\begin{proposition} \label{prop_laminar}
Let $T$ be a woven positive closed $(p,p)$-current on $V$. Let $K$ be a compact subset of $V$ such that $T$ has no mass on $K$. Assume that outside $K$ the current $T$ can be locally written as a wedge-product of positive closed $(1,1)$-currents with continuous potentials. If $2p\leq l$, we assume moreover that $T\wedge T=0$ on $X\setminus K$. Then $T$ is completely weakly laminar.
\end{proposition}
\proof
Observe that the case $2p>l$ can be reduced to the case $2p\leq l$ by replacing $V$ by $V\times V$ and $T$ by $T\otimes[V]$ in $V\times V$. Assume that $2p\leq l$.

Let $\nu$ be a measurable web associated with $T$. We have to show that $T$ is 
weakly laminar with respect to $\nu$.
We can refine $\nu$ in order to assume that $\nu$-almost every lame 
is defined by an irreducible  manifold which does not intersect $K$. Assume that $T$ is not completely weakly laminar and set $r:=l-p$. Then for a suitable $\nu$, there is a subset ${\cal W}$ of $\Lam_r^*(V)\times\Lam_r^*(V)$ of positive $\nu\times\nu$ measure such that for every $([Y],[Z])$ in ${\cal W}$ we have $Y\cap Z\not=\varnothing$ and  $\dim Y\cap Z=s$ for some integer $0\leq s\leq r-1$.  We can refine $\nu$ and reduce ${\cal W}$ in order to assume that all these sets $Y$ and $Z$ are closed submanifolds of a fixed open subset $\U$ of $V\setminus K$, as in the above local description of woven currents. 

The tangent cone of $Y\times Z$ with respect to the diagonal $\Delta$ of $V\times V$ is a non-empty variety. We deduce from the definition of tangent currents that no tangent current of $T\otimes T$ along $\Delta$ vanishes over $\U$. Here we identify $\U$ with an open subset of $\Delta$. Recall that such a tangent current is a positive closed current on the normal bundle to $\Delta$. It can be obtained locally as a limit value of the images of $T\otimes T$ by a sequence of dilations in the normal directions to $\Delta$, see \cite{DS12}.  
On the other hand, by Theorem 5.10 in \cite{DS12}, over $\U$, this tangent cone should be the pull-back of the current $T\wedge T$ to the normal vector bundle to $\Delta$. This contradicts the  hypothesis that $T\wedge T=0$. The proposition follows.
\endproof

Let $\Gr(V,r)$ denote the Grassmannian bundle over $V$ which is the set of points $(x,[v])$, where $x$ is a point in $V$ and $[v]$ is the direction of a simple complex tangent $r$-vector of $V$ at $x$. 
If $S$ is a current in $\Lam_r(V)$, write $S=\sum [\Lambda_i]$. We can lift each $\Lambda_i$ to $\Gr(V,r)$ by considering the set $\widehat \Lambda_i$ of points $(x,[v])$ with $x$ a regular point in $\Lambda_i$ and $v$ tangent to $\Lambda_i$ at $x$. If $\sum \|\widehat \Lambda_i\|$ is finite, $\widehat S:=\sum [\widehat \Lambda_i]$ is a current in $\Lam_r(\Gr(V,r))$ and we say that $\widehat S$ is the {\it lift} of $S$ to $\Gr(V,r)$. It does not depend on the choice of the decomposition $S=\sum [\Lambda_i]$.

Let $T$ be a woven positive closed current of bi-dimension $(r,r)$ and let $\nu$ be a measurable $r$-web associated with $T$. We have 
$$T=\int_{\Lam_r^*(V)} S d\nu(S).$$
Assume that $\nu$-almost every $S$ admits a lift to $\Gr(V,r)$. We can always have this property by refining $\nu$. 
If the integral 
$\int_{\Lam_r^*(V)} \|\widehat S\| d\nu(S)$ is finite, the current 
$$\widehat T:= \int_{\Lam_r^*(V)} \widehat S d\nu(S)$$
is well-defined and is called {\it a lift of $T$} to $\Gr(V,r)$. It may depend on the choice of  the measurable web $\nu$. The push-forward of $\widehat T$ to $V$ is always equal to $T$. 

\begin{definition} \rm \label{def_tame}
We say that $T$ is {\it almost tame} if it admits a measurable web $\nu$ as above with
$\int_{\Lam_r^*(V)} \|\widehat S\| d\nu(S)$ finite and if there is a positive {\bf closed} current $\widehat T'$ on $\Gr(V,r)$ such that 
$\widehat T'\geq \widehat T$ and the push-forward of $\widehat T'$ to $V$ is equal to $T$. 
We say that $T$ is {\it tame} if we can choose $\widehat T'$ equal to  $\widehat T$, i.e. the last current is closed. 
The above measurable web $\nu$ is said to be {\it almost tame} or {\it tame} respectively. 
\end{definition}

Note that almost tame currents are necessarily closed.

\begin{example} \label{ex_Fubini_bis}\rm
Consider the situation in Example \ref{ex_Fubini}. 
Set $\Omega_{k-r}:=\omega_\FS^{k-r}$.
 If $\widehat H$ is the lift of $H$ to $\Gr(\P^k,r)$, since the action of $\U(k+1)$ extends canonically to $\Gr(\P^k,r)$, the positive closed current
$$\widehat\Omega_{k-r}:=d^{-1}\int_{\tau\in\U(k+1)} \tau_*[\widehat H]d\sigma(\tau)$$
is a lift of $\Omega_{k-r}$ to $\Gr(\P^k,r)$. 
It is invariant under the action of $\U(k+1)$. 
Since $\Gr(\P^k,r)$ is a homogeneous space which is the quotient of $\U(k+1)$ by a subgroup, $\widehat\Omega_{k-r}$ is a smooth form. For $d=1$, we call $\widehat\Omega_{k-r}$ the {\it standard lift} of $\Omega_{k-r}$ to $\Gr(\P^k,r)$.
Note that if $H$ is a smooth hypersurface of degree $d$ of $\P^k$, then the lift $\widehat H$ of $H$ to $\Gr(\P^k,k-1)$ has volume of order $d^2$. This is the reason why the limit currents of varieties are not woven in general.   
\end{example}

Let $\pi:V\to W$ be a holomorphic submersion onto a compact complex manifold $W$. If $S=\sum [\Lambda_i]$ is as above write $S=S_1+S_2$ with $S_1:=\sum_1 [\Lambda_i]$ and $S_2:=\sum_2 [\Lambda_i]$, where $\sum_1$ is taken over  the $\Lambda_i$'s such that the restriction of $\pi$ to $\Lambda_i$ is generically of maximal rank and $\sum_2$ is taken over the other $\Lambda_i$'s. Let $T=\int_{\Lam^*_r(V)} [S] d\nu(S)$ be a woven positive closed current of bi-dimension $(r,r)$ associated with a measurable web $\nu$. 
Write $T=T_1+T_2$ with $T_i:=\int_{\Lam_r^*(V)} [S_i] d\nu(S)$. We have the following lemma that can be extended to the case where $\pi$ is a dominant meromorphic map. 

\begin{lemma} \label{lemma_tame_decom}
Let $T,\nu,T_1,T_2$ be as above. If $\nu$ is almost tame,  then $T_1$ and $T_2$ are closed.
\end{lemma}
\proof
We use the notations introduced above.
Denote by $Z$  the analytic set of points $(x,[v])$ in $\Gr(V,r)$  such that 
$v$ is not transverse to the fiber of $\pi$ through $x$.
Let $\widehat T_2$ be the restriction of $\widehat T'$ to $Z$. 
This is a positive closed current. 
Note that for every lame $S=\sum [\Lambda_i]$, we have $\widehat\Lambda_i \subset Z$ if and only if the rank of $\pi$ on $\Lambda_i$ is not maximal. 
Define $\widehat T_1:=\widehat T'-\widehat T_2$. This current is also positive and closed. By Definition \ref{def_tame}, the push-forwards of $\widehat T'$ and $\widehat T$  to $V$ are both equal to $T$. Therefore, $T_i$ is the push-forward of $\widehat T_i$ to $V$. The lemma follows. 
\endproof

\noindent
{\bf Woven currents as limits of analytic sets.} 
In dynamics, woven currents are often constructed as limits of currents of integration on analytic sets. 
The following result was obtained in \cite{Dinh2}.

\begin{theorem} \label{th_woven_limit}
Let $\Gamma_n$ be a sequence of analytic subsets of pure dimension $r$ in a projective manifold $V$ and let $d_n$ be positive numbers such $d_n^{-1}[\Gamma_n]$ converge to a current $T$. Let $\widehat \Gamma_n$ be the lift of $\Gamma_n$ to $\Gr(V,r)$. Assume that the $2r$-dimensional volume of $\widehat \Gamma_n$ is bounded by $cd_n$ for some constant $c>0$. Then $T$ is woven.
\end{theorem}

Note that we can lift the regular part of $\Gamma_n$ to $\Gr(V,r)$ and its compactification is an analytic subset of $\Gr(V,r)$ that we still denote by $\widehat \Gamma_n$. 

\medskip\noindent
{\it Sketch of the proof.} Since $V$ is projective, it can be embedded in a projective space. For simplicity, we can assume that $V$ is the projective space $\P^k$ and $d_n$ is the degree of $\Gamma_n$. Fix a generic central projection $\pi:\P^k\setminus I\to \P^r$, where $I$ is a projective subspace of dimension $k-r-1$ of $\P^k$ and $\P^r$ is identified with a projective subspace in $\P^k\setminus I$. If $z$ is a point in $\P^k\setminus I$, then $\pi(z)$ is the intersection of $\P^r$ with the projective subspace of dimension $k-r$ containing $I$ and $z$. 

If $z_0$ is a generic point in $\P^r$ and $U$ is a small neighbourhood of $z_0$, we can show that $\pi^{-1}(U)\cap \Gamma_n$ contains almost $d_n$ graphs over $U$ for $n$ large enough ($d_n$ is the maximal number one can have). This is the consequence of the property that the set of ramification of $\pi$ restricted to $\Gamma_n$ is small enough over $U$. We will see in Propositions \ref{prop_branch} and \ref{prop_branch_bis} below similar situations. The control of the ramification is obtained from the hypothesis on the volume $\widehat\Gamma_n$ using Fubini theorem and a generic choice of $\pi$, $z_0$. 

The limits of the obtained graphs as $n\to\infty$ form a part of $T$. We have to cover $\P^r$ with such open sets $U$ with different sizes in order to construct a complete measurable web associated with $T$. For the details, see  \cite{Dinh2}.
\hfill $\square$

\medskip

The following result can be deduced from the proof of the above theorem. 

\begin{proposition}\label{prop_woven_limit}
There is an increasing  sequence of integers $(n_i)$ and  measurable $r$-webs $\nu_{n_i}$ and $\nu$ associated with $d_{n_i}^{-1}[\Gamma_{n_i}]$ and $T$ such that $\nu_{n_i}\to\nu$ in the weak sense of measures on $\Lam_r^*(V)$. Moreover, we can write 
$\nu_{n_i}=\sum_{s\geq 1} \nu_{n_i}[s]$ and $\nu=\sum_{s\geq 1} \nu[s]$ such that $\nu_{n_i}[s]$ and $\nu[s]$ are $s$-elementary webs and $\nu_{n_i}[s]\to\nu[s]$ in the weak sense of measures on $\Lam^*_r(V)$.   
\end{proposition}

\begin{proposition} \label{prop_tame_limit}
Let $\Gamma_n, \widehat \Gamma_n$ and $T$  be as in the last theorem.  Assume that the $2r$-dimensional volume of the lift of $\widehat \Gamma_n$ to $\Gr(\Gr(V,r),r)$
 is bounded by $cd_n$ for some constant $c>0$. Then $T$ is almost tame. 
\end{proposition}
\proof
Let $S$ be a cluster value of the sequence $d_n^{-1}[\widehat \Gamma_n]$. Applying Theorem \ref{th_woven_limit} to $\widehat \Gamma_n$ implies that $S$ is woven. Here, in order to check the hypotheses of that theorem,  we need to lift $\Gamma_n$ twice.
With the construction explained above, we see that the elementary lames of $S$ are obtained as limits of open subsets of $\widehat \Gamma_n$. So the elementary lames whose projections on $V$ are of dimension $r$ are the lifts of some varieties in $V$ to $\Gr(V,r)$. For the other elementary lames,
their projections on $V$ vanish in the sense of currents.
It follows that for a suitable measurable web, the lift of $T$ to $\Gr(V,r)$ is bounded by $S$. Since $S$ is closed, the current $T$ is almost tame.
\endproof

Theorem \ref{th_woven_limit} can be extended to some local situation. We will need some steps in the proof of such local version that we recall below.

Let $\B_r$ denote the unit ball  and $\rho\B_r$ the ball of center 0 and of radius $\rho$ in $\C^r$. Consider an analytic subset 
$\Gamma$ of pure dimension $r$ of $3\B_r\times 3\B_s$, not necessarily irreducible, which is contained in $3\B_r\times 2\B_s$. For simplicity assume that $\Gamma$ is smooth.
So the natural projection from $\Gamma$ onto $3\B_r$ defines a ramified covering and we denote by $d$ its degree. 
The ramified locus is a divisor of $\Gamma$ with integer coefficients. Its push-forward to $3\B_r$ is a divisor with integer coefficients on $3\B_r$. The positive closed $(1,1)$-current associated with this divisor  is denoted by $[P]$ and is called the {\it postcritical current}.
Observe that when $[P]=0$ the set $\Gamma$ is a union of $d$ graphs over $3\B_r$.
The following result gives us a more quantitative property.

\begin{proposition} \label{prop_branch}
There is a constant $c>0$ independent of $\Gamma$ and $d$ such that $\Gamma\cap (2\B_r\times 3\B_s)$ contains at least 
$d-c\|P\|$ graphs over $2\B_r$.
\end{proposition}

The case of dimension $r=1$ is just a consequence of Riemann-Hurwitz's formula. The general case is reduced to the dimension 1 case by slicing $3\B_r$ by lines through the origin. We then stick graphs of dimension 1 in order to get graphs of dimension $r$ over $2\B_r$ using the main result in \cite{SW}, see \cite{Dinh2} for details.

Note that the proposition still holds for $\Gamma$ singular but the postcritical current $[P]$ has to be defined differently. 
The current $[P]$ depends strongly on the coordinate system we use. Therefore, in the general dynamical setting, we need a more subtle version of the last proposition. 

Define $\U:=4\B_r\times 3\B_s$ and assume now that $\Gamma$ is a smooth analytic subset of pure dimension $r$
of $\U$ which is contained in $4\B_r\times\B_s$ and is a ramified covering of degree $d$ over $4\B_r$.
Let $\Gr(\U,s)$ denote the set of point $(z,[v])$ where $z$ is a point in $\U$ and $[v]$ is the direction of a simple complex tangent  $s$-vector $v$ of $\U$ at $z$. It can be identified with the product of $\U$ with the Grassmannian $\G$ parametrizing the family of complex linear subspaces of dimension $s$ in $\C^{r+s}$ through a fixed point. 
We have $\dim\G=rs$ and $\dim \Gr(\U,s)=rs+r+s$.

We are interested in linear subspaces close enough to the vertical ones. More precisely, we consider linear subspaces parallel to a space of equation
$$z'=Az'' \quad \text{with} \quad z=(z',z'')\in \C^r\times\C^s,$$
where $A$ is a complex $r\times s$-matrix whose coefficients have modulus smaller than 1. The family of those matrices is identified with an open set $\G^\star$ in $\G$. Define $\Gr(\U,s)^\star:=\U\times \G^\star$. 

Denote by $\widetilde\Gamma$ the set of points $(z,[v])\in\Gr(\U,s)$ such that $z\in\Gamma$ and $v$ is not transverse to the tangent space of $\Gamma$ at $z$. So $\widetilde\Gamma$ is an analytic subset of $\Gr(\U,s)$ of pure dimension $rs+r-1$.
Indeed, it is not difficult to see that the restriction of $\widetilde\Gamma$ to a fiber $\{z\}\times\G$ is a hypersurface of this fiber. 
Define $\widetilde\Gamma^\star:=\widetilde\Gamma\cap \Gr(\U,s)^\star$. 

\begin{proposition} \label{prop_branch_bis}
There is a constant $c>0$ independent of $\Gamma$ and $d$ such that $\Gamma\cap (\B_r\times 3\B_s)$ contains at least 
$d-c\|\widetilde\Gamma^\star\|$ graphs over $\B_r$.
\end{proposition}

Fix a constant $\delta>0$ small enough depending only on $r$ and $s$. We will consider the projections 
$\pi_A:\C^r\times\C^s\to\C^r$ given by $\pi_A(z):=z'-Az''$ with $\|A\|\leq\delta$. We will apply Proposition \ref{prop_branch} to the coordinate system 
$$z_A=(z'_A,z''_A):=(z'-Az'',z'')$$ 
instead of $(z',z'')$. We will add the letter $A$ in the above notations when we use these coordinates.
In these coordinates, $\pi_A$ is just the natural projection on the first $r$ coordinates. Since $\delta$ is small, the following lemma is clear by continuity.

\begin{lemma}
The restriction of $\Gamma$ to $3\B_r^A\times 3\B_s^A$ is a ramified covering of degree $d$ over $3\B_r^A$ and is contained in $3\B_r^A\times  2\B_s^A$. If $Z$ is a graph over $2\B_r^A$ contained in $2\B_r^A\times 2\B_s^A$ then its restriction to $\B_r\times 3\B_s$ is also a graph over $\B_r$.  
\end{lemma}

\smallskip

\noindent
{\bf End of the proof of Proposition \ref{prop_branch_bis}.} 
Using the last lemma and Proposition \ref{prop_branch}, we only need to check that there is a matrix $A$ with $\|A\|\leq\delta$ such that the mass of the associated poscritical current $[P_A]$ is smaller than a constant times the mass of $\widetilde\Gamma^\star$. 
By Fubini's theorem, there is a matrix $A$ with $\|A\|\leq\delta$ such that the mass of $\widetilde\Gamma\cap (\U\times\{A\})$ is smaller than a constant times $\|\widetilde\Gamma^\star\|$, where the points in the last intersection are counted with multiplicity. It is enough now to observe that $P_A$ is equal on $3\B_r^A$ to the image of $\widetilde\Gamma\cap (\U\times\{A\})$ by $\pi_A$. This completes the proof of the proposition.
\hfill $\square$

\bigskip\noindent
{\bf Slicing theory for  woven currents.} In Section \ref{section_coh}, we already discussed slicing for positive closed currents and for varieties. Since woven currents are generated by pieces of varieties, the theory extends without difficulty to them. 
 
Let $\pi:V\to W$ and $l,m$ be as in the beginning of Section \ref{section_coh}. Assume for simplicity that $V$ and $W$ are projective manifolds. Let $T$ be a woven current on $V$ of bi-dimension $(r,r)$ associated with a web $\nu$. It follows from Federer's slicing theory for flat currents that  for almost every $y\in W$ the slice $\langle T|\pi|y\rangle$ exists and we have
$$\langle T|\pi|y\rangle =\int_{\Lam_r^*(V)} \langle S|\pi|y\rangle d\nu(S).$$
Note that the family of $y$ satisfying the above property depends not only on $T$ but also on the choice of $\nu$. 
When the above identity holds for $y$, we say that the web $\nu$ is {\it compatible} with the slice $\langle T|\pi|y\rangle$. A necessary condition for $\nu$ to be compatible with $\langle T|\pi|y\rangle$ is that $\nu$-almost every lame $S$  either is  disjoint from $\pi^{-1}(y)$ or intersects $\pi^{-1}(y)$ transversally at almost every point of intersection, see also Section \ref{section_coh}.
We describe now the situation that will be used later.

Let $\Gamma_n, T$ and  $d_n$  be as in 
Theorem \ref{th_woven_limit}.
Choose the covering of $V$ by cubes $\U_1,\ldots,\U_N$  as in Section 2.

\begin{proposition} \label{prop_slice_woven}
There is an increasing sequence of integers $(n_i)$ such that for almost every $y\in W$ the slices $\langle T|\pi|y\rangle$ and $d_{n_i}^{-1}\langle [\Gamma_{n_i}]|\pi|y\rangle$ are well-defined and $d_{n_i}^{-1}\langle [\Gamma_{n_i}]|\pi|y\rangle\to \langle T|\pi|y\rangle$. For almost every $y\in W$ there are measurable webs $\nu_{n_i}$, $\nu_{n_i}[s]$, $\nu$, $\nu[s]$, depending on $y$, which satisfy Proposition \ref{prop_woven_limit} and such that
 $\nu_{n_i}$ and $\nu$ are compatible with the above slices. Moreover,  for all $i,s$ and for every  graph $\Lambda$ corresponding to a point in the support of  $\nu_{n_i}[s]$ or  $\nu[s]$, the fat extension of $\Lambda$  either  is disjoint from  $\pi^{-1}(y)$ or intersects $\pi^{-1}(y)$ transversally. 
\end{proposition}
\proof
The first assertion is given by Proposition \ref{prop_slice_C0}. By replacing $(n_i)$ with a suitable subsequence and using Proposition \ref{prop_woven_limit}, we obtain  $\nu_{n_i}$, $\nu_{n_i}[s]$, $\nu$, $\nu[s]$
satisfying this proposition. Then for almost every $y$, the webs $\nu$ and $\nu_{n_i}$ are compatible with $\langle T|\pi|y\rangle$ and $\langle [\Gamma_{n_i}]|\pi|y\rangle$.   We explain how to modify the above webs in order to get the last property in the proposition. The modification depends on $y$. 

Recall that {\it almost} every lame used here  either is  disjoint from $\pi^{-1}(y)$ or intersects $\pi^{-1}(y)$ transversally at {\it almost} every point of intersection. We want in particular to remove the second word "almost" in the last sentence. For this purpose, we can simply remove from each lame a suitable proper analytic subset and this does not modify the currents. However, since we want the lames to be elementary, the problem is slightly more delicate. 

Observe that any elementary lame of order $s$ can be divided into elementary lames of order $s+1$ and gives us a refinement. 
Fix a constant $\epsilon_1>0$ small enough and consider the set $K$ of elementary lames $\Lambda$ of order 1 such that either the distance between the very fat extension of $\Lambda$ and $\pi^{-1}(y)$ is larger than $\epsilon_1$ or this very fat extension  intersects $\pi^{-1}(y)$ transversally and the angle between them is at least equal to $\epsilon_1$. This is a compact subset of $\Lam^*_r(V)$. 

We replace $\nu{[1]}$ with its restriction $\nu[1]_{|K}$ to $K$. We also refine $\nu[1]-\nu[1]_{|K}$ into a 2-elementary web that we add to $\nu[2]$. 
It is not difficult to find positive measures $\nu_{n_i}[1]'\leq \nu_{n_i}[1]$ such that $\nu_{n_i}[1]'\to \nu[1]_{|K}$
and for every $\Lambda$ in the support of $\nu_{n_i}'$  the fat extension of $\Lambda$ either is disjoint from $\pi^{-1}(y)$ or intersects $\pi^{-1}(y)$ transversally. We do a similar modification to $\nu_{n_i}[1]$ and $\nu_{n_i}[2]$: we replace $\nu_{n_i}[1]$ with $\nu_{n_i}[1]'$ and add to $\nu_{n_i}[2]$ the refinement of order 2 of $\nu_{n_i}[1]-\nu_{n_i}[1]'$. This completes the construction of $\nu[1]$ and $\nu_{n_i}[1]$.  

We repeat now the same modification for $\nu[2]$ and $\nu_{n_i}[2]$ using another $\epsilon_2$ small enough and the refinement of lames of order 2 into lames of order 3. We obtain the final webs $\nu[2]$ and $\nu_{n_i}[2]$. A simple induction allows to obtain the webs satisfying the proposition.  If we choose in each step the constant $\epsilon_s$ small enough, it is not difficult to insure that there is no mass lost after an infinite number of steps, that is, the currents associated with $\nu$, $\nu_{n_i}$, after modifications, are still equal to $T$, $d_{n_i}^{-1}[\Gamma_{n_i}]$ and compatible with the slicing by $\pi^{-1}(y)$. This ends the proof of the proposition.  
\endproof

We come back now to the situation in Section \ref{section_coh} where $V=V'\times \P^k$. We use the notations introduced there. Let $T$ be a woven current on $V$ of bi-dimension $(r,r)$ associated with a web $\nu$. Then for almost every 
$\xi\in \G$ the intersection $T\wedge [V'\times H_\xi]$ exists and we have 
$$T\wedge [V'\times H_\xi]=\int_{\Lam_r^*(V)} S\wedge [V'\times H_\xi]d\nu(S).$$
The family of $\xi$ satisfying the last identity depends on the choice of $\nu$.
When the above identity holds for $\xi$, we say that the lamination $\nu$ is {\it compatible} with the intersection $T\wedge [V'\times H_\xi]$. A necessary condition for $\nu$ to be compatible with $T\wedge [V'\times H_\xi]$ is that $\nu$-almost every lame $S$ either is  disjoint from $V'\times H_\xi$ or intersects $V'\times H_\xi$ transversally at almost every point of intersection.
We have the following direct consequence of Proposition \ref{prop_slice_C0} and Proposition \ref{prop_slice_woven}.

\begin{corollary} \label{cor_inter_woven}
Let $T$ and $\Gamma_n$ be as in Proposition \ref{prop_slice_woven}.
Then there is an increasing sequence of integers $(n_i)$ such that for almost every $\xi\in \G$ the intersections $T\wedge [V'\times H_\xi]$ and $d_{n_i}^{-1}[\Gamma_{n_i}]\wedge [V'\times H_\xi]$ are well-defined and 
we have $d_{n_i}^{-1}[\Gamma_{n_i}]\wedge [V'\times H_\xi]\to T\wedge [V'\times H_\xi]$. For almost every $\xi\in\G$, there are measurable webs $\nu_{n_i}$, $\nu_{n_i}[s]$, $\nu$, $\nu[s]$, depending on $\xi$, which satisfy 
Proposition \ref{prop_woven_limit}  and such that
$\nu_{n_i}$, $\nu$ are  compatible with the above intersections. Moreover, for all $i,s$ and for every  graph $\Lambda$ corresponding to a point in the support of  $\nu_{n_i}[s]$ or  $\nu[s]$,  the fat extension of $\Lambda$ either is disjoint from  $V'\times H_\xi$ or intersects $V'\times H_\xi$ transversally. 
\end{corollary}

We will also need the following lemma.
 
\begin{lemma} \label{lemma_cv_inter}
Let $\nu_n, \nu$ be $s$-elementary webs and let $T_n,T$ be the associated woven currents. Assume that  $\nu_n\to\nu$. Let $H$ be a projective subspace of  $\P^k$. Assume that for $\nu$-almost  every 
graph $\Lambda$, the fat extension of $\Lambda$ either is disjoint from $V'\times H$ or intersects $V'\times H$ transversally. Assume also that the intersections  $T_n\wedge [V'\times H]$ and  $T\wedge [V'\times H]$ exist and are compatible with the above webs. 
Then any limit value of $T_n\wedge [V'\times H]$ is larger than or equal to $T\wedge [V'\times H]$. 
\end{lemma}
\proof
We prove the lemma assuming only that $\lim \nu_n\geq \nu$. This allows to replace $\nu$ with its restrictions to suitable compact sets in order to get the above condition on $\Lambda$ for every graph $\Lambda$ corresponding to a point in the support of $\nu$. 
Let $\nu_n'$ be the restriction of $\nu_n$ to a small enough neighbourhood of the support of $\nu$. We have $\lim \nu_n'\geq \nu$ and $\nu_n'\leq \nu_n$. 
Observe that if a graph $\Lambda'$ is close enough to the graph $\Lambda$ in the lemma then the intersection $\Lambda'\cap (V'\times H)$ is transverse and close to $\Lambda\cap (V'\times H)$. It follows that for $n$ large enough, the current obtained by intersecting the lames of $\nu_n'$ with 
$[V'\times H]$ is close to $T\wedge [V'\times H]$. Since $\nu_n'\leq \nu_n$ and $\lim \nu_n'\geq \nu$ the lemma follows.
\endproof

\section{Green currents of H\'enon maps} \label{section_Green}

In this section, we give some crucial properties of Green currents associated with  H\'enon maps that we need for the proof of the main theorem. Let $f$ be a H\'enon map on $\C^k\subset\P^k$ as in the introduction.  
Define $T_+:=\tau_+^p$ and $T_-:=\tau_-^{k-p}$. They are positive closed currents of mass 1  respectively
 of bi-degree $(p,p)$ and $(k-p,k-p)$ with support in $\overline K_+$ and $\overline K_-$. We have $f^*(T_+)=dT_+$ and $f_*(T_-)=d T_-$. 
For dimension reason $T_\pm$ have no mass on $I_\pm$ and  then they have no mass on the hyperplane at infinity $H_\infty$ because  $\overline K_\pm=K_\pm\cup I_\pm$.
The following result was obtained in \cite{DS10} using the theory of super-potentials.

\begin{theorem} \label{th_unique_Green}
The current $T_+$ is the unique positive closed $(p,p)$-current of mass $1$ supported by $\overline K_+$.  In particular, it is an extremal positive closed $(p,p)$-current on $\P^k$.
Moreover, if $S$ is a positive closed $(p,p)$-current of mass $1$, smooth in a neighbourhood of $I_-$, then $d^{-n}(f^n)^*(S)$ converges to $T_+$ as $n\to\infty$. 
\end{theorem}

Note  that
a similar version of the above theorem holds for $T_-,\overline K_-$ and $f^{-1}$. 
The following theorem generalizes some results  in \cite{BLS} for dimension $k=2$ and  \cite{Dinh2} for higher dimension.

\begin{theorem} \label{th_lam_Green}
The currents $T_+$ and $T_-$ are completely weakly laminar and tame.
\end{theorem}
\proof We will prove the result for $T_+$. The case of $T_-$ can be obtained in the same way.
Let $L$ be a linear subspace of dimension $k-p$ of $\P^k$ which does not intersect $I_-$. It follows from Theorem \ref{th_unique_Green}  that $d^{-n} [f^{-n}(L)]$ converges to $T_+$. Denote by $L_n$ the compactification of $f^{-n}(L)\cap \C^k$ and let $\widehat L,\widehat L_n$ be respectively the lifts of $L,L_n$ to $\Gr(\P^k,k-p)$, see the notations in Section \ref{section_woven} and the end of Section \ref{section_coh}.
We have $\widehat L_n=\widehat f^{-n}(\widehat L)$ in the Zariski open subset $\C^k\times\G$ of $\Gr(\P^k,k-p)$. The varieties $\widehat L$ and $\widehat L_n$ are irreducible and not contained in the complement of $\C^k\times\G$. 

The properties of $f^n$ described in Example \ref{ex_Henon} imply that the volume of $\widehat L_n$ is equal to $O(d^n)$. We used here that the mass of a positive closed current, in particular the volume of an analytic set, only depends on its cohomology class. 
In the same example, we see that the volume of the lift of $\widehat L_n$ to $\Gr(\Gr(\P^k,k-p),k-p)$ is also equal to $O(d^n)$. Therefore, by Theorem \ref{th_woven_limit} and Proposition \ref{prop_tame_limit}, the Green current $T_+$ and the cluster values of $d^{-n}[\widehat L_n]$ are almost tame woven positive closed currents. 

Let $S$ be such a cluster value on $\C^k\times \G$. Since it has finite mass, its extension by 0 is a positive closed current on $\Gr(\P^k,k-p)$.  
Recall that we can construct a measurable web associated with $S$ such that its lames are elementary and  obtained as limits of open subsets of $\widehat L_n$. Write $S=S'+S''$, where  $S'$ (resp. $S''$) is formed  by lames whose projections on $\C^k$ are non-degenerate (resp. degenerate), i.e. of dimension $k-p$ (resp. smaller than $k-p$).  
So for a suitable measurable web, $S'$ is a lift of $T_+$ to $\Gr(\P^k,k-p)$. By Lemma \ref{lemma_tame_decom}, $S'$ and $S''$ are closed.

We want to prove that $T_+$ is tame. For this purpose, it is sufficient to check that $S''=0$.  Consider the family $\Fc$ of all currents $S''$ obtained as above for different cluster values $S$. These currents have a bounded mass. If $S$ is the limit of a sequence $d^{-n_j} [\widehat L_{n_j}]$ and $S_1$ is a cluster value of   $d^{-n_j+1} [\widehat L_{n_j-1}]$ then $d^{-1}\widehat f^\bullet (S_1)=S$ on $\C^k\times\G$ and hence on $\Gr(\P^k,k-p)$. It follows that there is a current $S_1''$ in $\Fc$ such that $S''=d^{-1}\widehat f^\bullet (S_1'')$.  By induction, there are currents $S_n''$ in $\Fc$ such that $d^{-n}(\widehat f^n)^\bullet (S_n'')=S''$.
The h-dimension of currents in $\Fc$ with respect to the projection on $\P^k$ is smaller than $k-p$. 
Therefore, their  cohomology classes are in $E_{q,p+1}$ for some $q$, see Section \ref{section_coh} for the notation. The claim in the proof of Proposition \ref{prop_coh_prod} implies that  the norm of $(\widehat f^n)^*$ on $E_{q,p+1}$ is equal to $o(d^n)$. The relation between $S''$ and $S_n''$ and the fact that $S_n''$ have a bounded mass imply that $S''=0$.  So $T_+$ is tame.

It remains to prove that $T_+$ is completely weakly laminar. Recall that $T_+=\tau_+^p$ and $\tau_+$ has local continuous potentials outside $I_+$. Moreover, $T_+$ has no mass on $I_+$ and $\tau_+^{p+1}=0$ outside $I_+$, see \cite{Sibony} for details. Proposition \ref{prop_laminar} implies the result. 
\endproof

Recall that  the measure $\mu$ is mixing and hyperbolic,  see \cite{deThelin2, Dinh}. It admits $p$ positive and $k-p$ negative Lyapounov exponents. For $\mu$ almost every point $z$ denote by $E_u(z)$ the unstable tangent subspace of $\P^k$ at the point $z$ and $E_s(z)$ the stable one.  So we have $\dim E_u(z)=p$ and $\dim E_s(z)=k-p$.
We will denote by $[E_s(z)]$ the direction of the complex tangent $(k-p)$-vectors at $z$ defining $E_s(z)$. 
The set of points $(z,[E_s(z)])$ in  $\Gr(\P^k,k-p)$ can be seen as a measurable graph over $\mu$-almost every point in the support of $\mu$. So we can lift $\mu$ to a probability measure $\mu_+$ on this graph that we call {\it stable Oseledec measure}  associated with $f$ and $\mu$. Since the stable bundle is invariant and $\mu$ is mixing, $\mu_+$ is also invariant under $\widehat f$ and mixing.  We can construct in the same way the {\it unstable Oseledec measure} $\mu_-$ associated with $f$ and $\mu$. It is a probability measure  on the set of points $(z,[E_u(z)])$ in $\Gr(\P^k,p)$.

The following result characterizes Oseledec measures. 
We will use it for $q=0$ and for the form $\widehat\Omega_p$ defined in Example \ref{ex_Fubini_bis} in order to prove the intersection properties we need. 
Recall that over $\C^k$ we identify $\Gr(\P^k,k-p)$ with $\C^k\times \G$. Let $m$ denote the dimension of $\G$. Let $\pi:\Gr(\P^k,k-p)\to\P^k$ denote the canonical projection.

\begin{proposition} \label{prop_main_inter}
Let $q$ be an integer such that $0\leq q\leq p$. Let $\alpha$ be a smooth positive closed form of bi-dimension $(k-p+q,k-p+q)$  on  $\Gr(\P^k,k-p)$.  Let $c$ be the mass of the current $\pi_*(\alpha)$.
Then 
$d_+^{-(p-q)n}(\widehat f^n)^* (\alpha)\wedge \pi^*(\tau_+^q\wedge T_-)$ defines a positive measure of mass $c$ on $\C^k\times \G$ which converges to $c\mu_+$ as $n\to\infty$.
\end{proposition}

Since $\pi$ is a submersion, $\pi_*(\alpha)$ is a smooth positive closed $(p-q,p-q)$-form. 
Since $(\widehat f^n)^* (\alpha)$ is smooth on  $\C^k\times\G$, the wedge-product in the proposition is well-defined and is a positive measure.  Its mass is equal to the mass of its push-forward to $\P^k$. This push-forward is equal to 
$d_+^{-(p-q)n}(f^n)^*  \pi_*(\alpha)\wedge \tau_+^q\wedge T_-$. 
Since $(f^n)^*  \pi_*(\alpha)$ is smooth on the support of $T_-$, the last wedge-product on $\C^k$ coincides with the same wedge-product on $\P^k$. 
Since $\tau_+$ and $T$ are of mass 1 and 
$d_+^{-(p-q)n}(f^n)^*$ is equal to the identity on $H^{p-q,p-q}(\P^k,\C)$, it is not difficult to see that 
the mass of the considered measure is equal to the mass $c$ of $\pi_*(\alpha)$. 

For the proof of the last proposition, we will use a decreasing induction on $q$. The following lemma proves the case $q=p$. The main ingredient here is  Oseledec's theorem.

\begin{lemma}\label{lemma_Oseledec}
Let $\alpha$ be a smooth  closed $(m,m)$-form on $\Gr(\P^k,k-p)$. Then we have $(\widehat f^n)^* (\alpha)\wedge \pi^*(\mu)\to c\mu_+$ on $\C^k\times\G$, where $c$ is the constant equal to $\pi_*(\alpha)$. 
\end{lemma}
\proof
Since $\alpha$ can be written as a combination of smooth positive closed forms, we can assume that it is positive. 
Observe that $\pi_*(\alpha)$ is a closed $(0,0)$-current. So it is defined by the above constant $c$. For simplicity, assume that $c=1$. We deduce that for any probability measure $\nu$ on $\C^k$ the positive measure $\alpha\wedge \pi^*(\nu)$ is of mass $1$ or equivalently the push-forward of $\alpha\wedge \pi^*(\nu)$ to $\C^k$ is equal to $\nu$. 
In particular, the mass of $\alpha\wedge [\pi^{-1}(z)]$ is $1$ for every $z\in \C^k$. 

Let $\widehat\varphi$ be a smooth test function on $\Gr(\P^k,k-p)$. We have to show that 
$$\big\langle (\widehat f^n)^\bullet (\alpha)\wedge \pi^*(\mu),\widehat\varphi\big\rangle \to \langle \mu_+,\widehat\varphi\rangle.$$ 
Define for $\mu$-almost every $z$ the function $\psi(z):=\widehat\varphi(z,[E_s(z)])$. Define also $\widehat\psi:=\psi\circ\pi$. 
For $\mu$-almost every point $z$, denote by $\Sigma(z)$ the set of points $(z,[v])$ in $\pi^{-1}(z)$ such that $v$ is not transverse to $E_u(z)$. It is a hypersurface  in $\pi^{-1}(z)$ and hence of Lebesgue measure 0 there. 
By Oseledec's theorem \cite[p.234]{Walters}, if $(z,[v])$ is out of $\Sigma(z)$ and $z_n:=f^{-n}(z)$, then the distance between $\widehat f^{-n}(z,[v])$ and $(z_n,[E_s(z_n)])$ tends to 0. It follows that $\widehat\varphi\circ \widehat f^{-n}-\widehat\psi\circ \widehat f^{-n}$ tends to 0 outside the union of $\Sigma(z)$. Since this function is bounded, using that $\alpha$ is smooth, we deduce that 
$$\big\langle \alpha\wedge \pi^*(\mu), \widehat\varphi\circ \widehat f^{-n}-\widehat\psi\circ \widehat f^{-n}\big\rangle \to 0.$$
This together with the following identities imply the result.

We have  since $\mu$ is invariant
$$\big\langle \alpha\wedge \pi^*(\mu), \widehat\varphi\circ \widehat f^{-n}\big\rangle= \big\langle (\widehat f^n)^\bullet\big(\alpha\wedge \pi^*(\mu)\big),\widehat\varphi\big\rangle
= \big\langle (\widehat f^n)^\bullet(\alpha)\wedge \pi^*(\mu),\widehat\varphi\big\rangle$$
and since the push-forward of $\alpha\wedge \pi^*(\mu)$ to $\C^k$ is equal to $\mu$
$$\big\langle \alpha\wedge \pi^*(\mu), \widehat\psi\circ \widehat f^{-n}\big\rangle
= \big\langle \alpha\wedge \pi^*(\mu), \pi^*(\psi\circ f^{-n})\big\rangle = \langle \mu, \psi\circ f^{-n}\rangle=\langle\mu,\psi\rangle.$$
The last integral is also equal to $\langle\mu_+,\widehat\varphi\rangle$. The proposition follows.
\endproof

\begin{lemma} \label{lemma_mu+}
Let $\widehat \mu$ be a positive measure on $\C^k\times \G$. 
Assume there is a negative current $U$ of bi-dimension $(1,1)$ on $\C^k\times \G$ such that $\widehat\mu-\mu_+=\ddc U$ and $\pi_*(U)=0$. Then $\widehat\mu=\mu_+$.  
\end{lemma}
\proof
By hypotheses, we have $\pi_*(\widehat\mu)=\pi_*(\mu_+)$. Therefore, $\widehat\mu$ is a probability measure and $\pi_*(\widehat \mu)=\pi_*(\mu_+)=\mu$. Thus, we can write 
$\widehat\mu=\int\widehat\mu_z d\mu(z)$ and $\mu_+=\int \mu_{+,z} d\mu(z)$, where $\widehat\mu_z$ and $\mu_{+,z}$ are respectively the conditional measures of $\widehat\mu$ and $\mu_+$ with respect to $\pi$ and $\mu$.
Note that  $\widehat\mu_z$ and $\mu_{+,z}$ are probability measures on $\pi^{-1}(z)$ which are defined for $\mu$-almost every $z$. It is enough to prove that  $\widehat\mu_z=\mu_{+,z}$. 

Since $\pi_*(U)=0$, as in Section 3 of \cite{DS12}, we obtain that $U$ is a vertical current with respect to $\pi$ in the sense that it 
 can be decomposed into currents on fibers of $\pi$. More precisely, there is a positive measure $\mu'$ on $\C^k$ and for $\mu'$-almost every point $z\in \C^k$ there is a negative current $U_z$ on $\pi^{-1}(z)$ such that $U=\int U_z d\mu'(z)$.  Write $\mu'=\mu'_r+\mu'_s$ with $\mu'_r$ absolutely continuous and $\mu'_s$ singular with respect to $\mu$. Multiplying $U_z$ by a suitable constant depending on $z$, possibly by 0, allows to assume that $\mu'_r=\mu$. 

\medskip\noindent
{\bf Claim.} We have $\ddc U_z=\widehat\mu_z-\mu_{+,z}$ for $\mu$-almost every $z$ and $\ddc U_z=0$ for $\mu'_s$-almost every $z$. 

\medskip

Assuming the claim, we first complete the proof of the lemma.
By definition, $\mu_{+,z}$ is a Dirac mass. For simplicity, we use a local coordinate system $x$ on $\pi^{-1}(z)$ so that $\mu_{+,z}$ is the Dirac mass $\delta_0$ at 0. 
We only have to check that $\widehat \mu_z\geq \delta_0$.
Choose a negative function $\varphi$ with support in a small neighbourhood of 0 which is smooth outside 0 and equal to $\log\|x\|$ near 0. Define $\varphi_n:=\max(\varphi,-n)$. 

Assume that the inequality $\widehat\mu_z\geq \delta_0$ does not hold. 
Since $\widehat\mu_z$ and $\mu_{+,z}$ have the same mass, we easily deduce that $\langle \widehat\mu_z-\mu_{+,z},\varphi_n\rangle\to\infty$ as $n\to\infty$. On the other hand, we have
$$\langle \widehat\mu_z-\mu_{+,z},\varphi_n\rangle =\langle \ddc U_z,\varphi_n\rangle =\langle U_z,\ddc \varphi_n\rangle.$$
The last integral is bounded above because $U_z$ is negative and $\ddc\varphi_n$ is positive in a fixed neighbourhood of 0 and smooth outside this neighbourhood. This is a contradiction. 

It remains to prove the claim. Define $V_z:=\ddc U_z-\widehat\mu_z+\mu_{+,z}$ for $\mu$-almost every $z$ and $V_z:=\ddc U_z$ for $\mu'_s$-almost every $z$. We have to show that $V_z=0$ for $\mu'$-almost every $z$.
Consider a dense sequence $\phi_n$ in the space of smooth test functions with compact supports in $\C^k\times\G$. 
It is enough to check for each $\phi:=\phi_n$ that $\langle V_z,\phi\rangle=0$ for $\mu'$-almost every $z$ or equivalently for any smooth function $\chi(z)$ with compact support in $\C^k$
$$\int \langle V_z,\phi\rangle \chi(z) d\mu'(z)=0.$$
By definition of $V_z$, the last integral is equal to $\big\langle \ddc U-\widehat \mu+\mu_+,(\chi\circ\pi)\phi \big\rangle$.
Therefore, it vanishes by hypotheses. This completes the proof of the lemma.
\endproof

\begin{lemma} \label{lemma_U_faible} 
Let $S$ be a positive closed $(q,q)$-current on $\Gr(\P^k,k-p)$. Let $\gamma$ be a   continuous negative form of bi-dimension $(k-p+q+1,k-p+q+1)$ on a neighbourhood of $\pi^{-1}(\overline K_-)$.   Then the mass of $d^{-n}(\widehat f^n)^\bullet(S\wedge \gamma)\wedge \pi^*(T_-)$ on $\C^k\times \G$ is bounded by a constant independent  of $n$.  
 \end{lemma}
\proof
Observe that the form $(\widehat f^n)^*(\gamma)$ is defined and continuous on a neighbourhood of $\pi^{-1}(K_-)$. 
Moreover, $T_-$ is a power of a positive closed $(1,1)$-current with continuous potential on $\C^k$. Therefore, the wedge-product in the lemma is well-defined on $\C^k\times\G$. 
Choose a smooth positive closed form $\theta$ on $\Gr(\P^k,k-p)$ such that $\gamma\geq -\theta$ on a neighbourhood of $\pi^{-1}(\overline K_-)$ which contains the support of $\pi^*(T_-)$. Replacing $\gamma$ with $-\theta$ allows to assume that $\gamma$ is negative closed and smooth on $\Gr(\P^k,k-p)$.  

It is now clear that the mass of $d^{-n}(\widehat f^n)^\bullet(S\wedge \gamma)$ is bounded since $d^{-n}(\widehat f^n)^*$ is bounded on the Hodge cohomology of $\Gr(\P^k,k-p)$, see Example \ref{ex_Henon}. It follows that the mass of $d^{-n}(\widehat f^n)^\bullet(S\wedge\gamma)\wedge \pi^*(T_-)$
is also bounded. To see this last point, it is enough 
to use that $T_-=\tau_-^{k-p}$ and to approximate $\tau_-$ by a sequence of smooth positive closed $(1,1)$-forms on $\P^k$ with decreasing potentials. Recall that the mass of a positive closed current on a compact K\"ahler manifold depends only on its cohomology class. 
The lemma follows.
\endproof

\noindent
{\bf End of the proof of Proposition \ref{prop_main_inter}.}
We use a decreasing induction on $q$. The case $q=p$ was considered in Lemma \ref{lemma_Oseledec}. 
Assume the proposition for $q+1$ with some $0\leq q\leq p-1$. We show it  for $q$. 
For simplicity assume that the mass $c$ of $\pi_*(\alpha)$ is 1. 
Observe that the bi-dimension of $\alpha$ is at most equal to $(k-1,k-1)$. Therefore, Leray's spectral theory applied to the fibration by $\pi$ gives a formula similar to the K\"unneth formula in the product case and implies that the cohomology class $\{\alpha\}$ 
belongs to the ideal generated by $\pi^*(\oplus_{r\geq 1}H^{r,r}(\P^k,\C))$, see \cite{Voisin}. Since the class $\{\tau_+\}$ and its powers generate the ideal 
$\oplus_{r\geq 1} H^{r,r}(\P^k,\C)$, we deduce that $\{\alpha\}$
can be written as the cup-product of 
$\{\pi^*(\tau_+)\}$ with the  class of some smooth closed form $\beta$ of the right bi-degree. We can assume that $\beta$ is positive since we can always reduce the problem to this case by considering linear combinations of positive forms. The condition $c=1$ implies that the mass of $\pi_*(\beta)$ is 1.

Since $\tau_+$ is of mass 1, it belongs to the class of the Fubini-Study form $\omega_\FS$. We deduce that $\alpha$ is cohomologous to $\pi^*(\omega_\FS)\wedge\beta$. In particular, there is a smooth form $\gamma'$ such that 
$\alpha=\pi^*(\omega_\FS)\wedge \beta+\ddc \gamma'$. Adding to $\gamma'$ a suitable negative closed form allows to assume that $\gamma'$ is negative. There is also a quasi-psh function $u$ such that $\omega_\FS=\tau_+-\ddc u$. Since $\tau_+$ has continuous potential outside $I_+$, the function $u$ is continuous outside $I_+$. Adding to $u$ a suitable constant allows to assume that it is positive on a neighbourhood of $\overline K_-$. Therefore, we have $\alpha=\pi^*(\tau_+)\wedge\beta+\ddc\gamma$ with $\gamma:=\gamma'-(u\circ\pi)\beta$.

Using that $d_+^{-1}f^*(\tau_+)=\tau_+$ and $\widehat f^\bullet\circ \pi^*=\pi^*\circ f^\bullet$, we obtain 
$$d_+^{-(p-q)n}(\widehat f^n)^\bullet (\alpha)\wedge \pi^*(\tau_+^q\wedge T_-) =  
d_+^{-(p-q-1)n}(\widehat f^n)^\bullet (\beta)\wedge \pi^*(\tau_+^{q+1}\wedge T_-) 
+ \ddc U_n$$
with $U_n:=d^{-n}(\widehat f^n)^\bullet\big(\pi^*(\tau_+^q)\wedge \gamma\big)\wedge \pi^*(T_-)$. 
The induction hypothesis implies that the first term in the last sum converges to $\mu_+$. By Lemma \ref{lemma_U_faible}, the negative currents $U_n$ have bounded masses.  So we can extract convergent subsequences.
We claim that the sequence of measures
$U_n\wedge \pi^*(\omega_\FS)$ converges to 0. This property implies that all cluster values $U$ of $U_n$ satisfy $\pi_*(U)=0$ and 
 Lemma \ref{lemma_mu+} gives the result.

It remains to prove the claim. As in the proof of Lemma \ref{lemma_U_faible}, we can assume that $\gamma$ is a negative closed smooth form. Since $\widehat f^n$ is an automorphism of $\C^k\times \G$, the mass of the measure $U_n\wedge \pi^*(\omega_\FS)$ is equal 
to the mass of its image by $\widehat f^n$, i.e. the mass of 
$$\pi^*(\tau_+^q)\wedge \gamma\wedge d^{-n}(\widehat f^n)_\bullet\pi^*(T_-\wedge \omega_\FS)=\pi^*(\tau_+^q)\wedge \gamma\wedge \pi^*\big[d^{-n}(f^n)_\bullet(T_-\wedge \omega_\FS)\big].$$
Since $(f^n)_*$ acts on $H^{k-p+1,k-p+1}(\P^k,\C)$ as the multiplication by $d_+^{(p-1)n}$, the current $d^{-n}(f^n)_\bullet(T_-\wedge \omega_\FS)$ tends to 0 as $n\to\infty$. The result follows.
\hfill $\square$

\medskip

\begin{remark} \rm
We don't know in general if $d^{-n}(\widehat f^n)^\bullet \big(\alpha\wedge \pi^*(\tau_+^q)\big)$ converges to a constant times a canonical lift of $T_+$ to $\Gr(\P^k,k-p)$ even for $q=0$ or $q=p$. We don't know if $T_+$ admits a unique lift to $\Gr(\P^k,k-p)$. This is true when $p=1$ because in this case we can show that $T_+$ is completely laminar. 
\end{remark}

\section{Equidistribution of saddle periodic points} \label{section_equi}

Let $F:\C^k\times\C^k\to \C^k\times\C^k$ be the polynomial automorphism defined in Example \ref{ex_Henon}. We extend it to a bi-rational map of $\P^k\times\P^k$. Recall that 
Proposition \ref{prop_ex} can be applied to $F$.
The indeterminacy sets of $F$ and $F^{-1}$ are  $(I_+\times \P^k)\cup (\P^k\times I_-)$  and $(I_-\times\P^k)\cup (\P^k\times I_+)$ respectively. Its dynamics is similar to the one of H\'enon maps. 
 Let $\pi_1,\pi_2:\P^k\times\P^k\to\P^k$ denote the canonical projections.
 The following result is proved in the same way as for  Theorem \ref{th_unique_Green}. We can also deduce it from that theorem in $\P^k$ and Propositions \ref{prop_current_tensor} and \ref{prop_product_map}. 

\begin{theorem} \label{th_F_equi}
Let $S$ be a positive closed $(k,k)$-current on $\P^k\times\P^k$. Assume that the support of $S$  does not intersect $I_-\times I_+$. Then 
$d^{-2n} (F^n)^\bullet(S)$ converges to $c\T_+$ as $n\to\infty$, where $\T_+:=T_+\otimes T_-$ and $c:=\big\langle S,\pi_1^*(\omega_\FS^{k-p})\wedge \pi_2^*(\omega_\FS^p)\big\rangle$. A similar result holds for $F^{-1}$. 
\end{theorem}

We will use the following corollary.

\begin{corollary}
Let $\Gamma_n$ denote the closure of the graph of $f^n$ in $\P^k\times \P^k$. Then the sequence of positive closed $(k,k)$-currents $d^{-n}[\Gamma_n]$ converges to $\T_+$ as $n\to\infty$.
\end{corollary}
\proof
If $n$ is even, we can write $d^{-n}[\Gamma_n]=d^{-n}(F^{n/2})^\bullet[\Delta]$. Otherwise, write $d^{-n}[\Gamma_n]=d^{-n+1}(F^{(n-1)/2})^\bullet(d^{-1}[\Gamma_1])$.  Since $f$ is a H\'enon map, it is not difficult to check that $\Delta$ and $\Gamma_n$ do not intersect $I_-\times I_+$ for $n\geq 0$. Therefore, we can apply Theorem \ref{th_F_equi}. Since $\Gamma_n$ is the closure of the graph of $f^n$, the constant $c$ there is equal to
$$d^{-n}\int_{\Gamma_n} \pi_1^*(\omega_\FS^{k-p})\wedge \pi_2^*(\omega_\FS^p)=d^{-n}\int_{\P^k}\omega_\FS^{k-p}\wedge (f^n)^*(\omega_\FS^p).$$
The last expression defines the mass of the current $d^{-n}(f^n)^*(\omega_\FS^p)$ which is equal to 1 because the action of $f^*$ on $H^{p,p}(\P^k,\C)$ is the multiplication by $d$. The corollary follows.
\endproof

Denote by $\widehat \Gamma_n$ and $\widehat \Delta$ the lifts of $\Gamma_n$ and  $\Delta$ to $\Gr(\P^k\times\P^k,k)$. 
Let $\widehat F$ denote the canonical lift of $F$ to $\Gr(\P^k\times\P^k,k)$. We have $\widehat\Gamma_n=\widehat F^{-n/2}(\widehat\Delta)$ if $n$ is even and $\widehat\Gamma_n=\widehat F^{-(n-1)/2}(\widehat\Gamma_1)$ if $n$ is odd.
Let $(n_i)$ be an increasing sequence of integers such that $d^{-n_i}[\widehat\Gamma_{n_i}]$ converges to a positive closed $(k,k)$-current $\widehat \T_+$. As in the proof of Theorem \ref{th_lam_Green}, we obtain the following result.

\begin{proposition} \label{prop_F_laminar}
The current $\T_+$ is completely weakly laminar and tame.
The  current $\widehat \T_+$ is a lift of $\T_+$ to $\Gr(\P^k\times\P^k,k)$. In particular, we have $\Pi_*(\widehat \T_+)=\T_+$, where 
$\Pi:\Gr(\P^k\times\P^k,k)\to\P^k\times\P^k$ denotes the canonical projection. 
\end{proposition}

Note that over $\C^k\times\C^k$, $\Gr(\P^k\times\P^k,k)$ can be identified with the product of $\C^k\times\C^k$ with a Grassmannian $\GG$. 
Recall that for $\mu$-almost every point $z$ we denoted by $E_u(z)$ the unstable tangent subspace and $E_s(z)$ the stable one. Denote by $\widehat\mu^\Delta$ the lift of $\mu$ to the set of points $(z,z,[E_s(z)\times E_u(z)])$ in $\Gr(\P^k\times\P^k,k)$. Here is a key point in the proof of our main result.

\begin{proposition} \label{prop_a_weak_intersection}
We have that $\widehat \T_+\curlywedge \Pi^*[\Delta]$ is well-defined and equal to $\widehat\mu^\Delta$.
\end{proposition}

The existence of the last  wedge-product means that 
\begin{enumerate}
\item[(a)] There is a unique tangent current of $\widehat \T_+$ along $\Pi^{-1}(\Delta)$. More precisely, when one dilates local coordinates in the normal directions to $\Pi^{-1}(\Delta)$, the image of $\widehat \T_+$ converges to a unique positive closed current on the normal vector bundle of $\Pi^{-1}(\Delta)$ independently of the choice of coordinates. The limit is called tangent current.
\item[(b)] The tangential h-dimension of $\widehat \T_+$ along $\Pi^{-1}(\Delta)$ is minimal, i.e. 0 in our case. Equivalently, in our case, the tangent current can be decomposed into currents of integration on fibers of the above normal vector bundle.
\end{enumerate}
In such a situation, the tangent current is the pull-back to the above normal vector bundle of a positive measure which is by definition the weak intersection $\widehat \T_+\curlywedge \Pi^*[\Delta]$.
We refer to \cite{DS12} for details.

The uniqueness of this tangent current requires the laminar properties of $\widehat\T_+$ and other dynamical arguments. We will analyze $\widehat\T_+$ using various projections.  
We first prove the property on the tangential h-dimension of $\widehat\T_+$. See \cite{DS12} for the density $\kappa_r$ of positive closed currents. It describes the cohomology classes of the tangent currents along a submanifold, according to their behavior along various directions. Our study in \cite{DS12} shows that tangent currents are not unique in general but they are in the same cohomology class of the normal vector bundle to the submanifold. We then obtain the cohomology classes $\kappa_r$ on the submanifold using Leray's theory.

\begin{lemma} \label{lemma_kappa_T+}
We have $\kappa_r(\widehat \T_+,\Pi^{-1}(\Delta))=0$ if $r>0$, i.e. the tangential h-dimension of $\widehat\T_+$ along $\Pi^{-1}(\Delta)$ is $0$.
\end{lemma}
\proof
Assume there is an integer $r\geq 1$ such that $\kappa_r(\widehat \T_+,\Pi^{-1}(\Delta))\not =0$. Choose $r$ maximal satisfying this property. Recall that such a maximal $r$ is called the tangential h-dimension of $\widehat\T_+$ along $\Pi^{-1}(\Delta)$.  
By  Lemma 3.8 in \cite{DS12}, $\kappa_r(\widehat \T_+,\Pi^{-1}(\Delta))$ is a pseudo-effective cohomology class of bi-dimension $(r,r)$ on $\Gr(\P^k\times\P^k,k)$ that can be represented by a positive closed current in $\supp(\T_+)\cap\Pi^{-1}(\Delta)$. 
Observe that the intersection $\supp(\T_+)\cap \Pi^{-1}(\Delta)$ has a compact projection in $\C^k\times\C^k$.
So we only need to work in $\C^k\times\C^k\times \GG$.

Define $\widehat \T_{+m}:=d^{-2m} (F^m)^\bullet (\widehat \T_+)$.  Since $\widehat F$ is an automorphism on $\C^k\times\C^k\times\GG$, we have
$\kappa_r(\widehat \T_+,\Pi^{-1}[\Delta])=\kappa_r(\widehat \T_{+m},d^{-2m}\Pi^*[\Gamma_{-2m}])$. Theorem \ref{th_F_equi} applied to $F^{-1}$ implies that $d^{-2m}[\Gamma_{-2m}]$ converges to $\T_-:=T_-\otimes T_+$. So $d^{-2m}\Pi^*[\Gamma_{-2m}]$ converges to $\Pi^*(\T_-)$. 
Let $\widehat \T_+'$ denote a limit value of the sequence
 $\widehat \T_{+m}$. Since $\Pi^*(\T_-)$ is the wedge-product of $(1,1)$-currents with continuous potentials, the intersection $\widehat \T_+'\wedge \Pi^*(\T_-)$ is a well-defined  positive measure. It follows from the theory of densities that the density dimension between $\widehat \T_+'$ and $ \Pi^*(\T_-)$ is zero. Corollary 5.8 in \cite{DS12} on the upper semi-continuity of densities implies that 
$\kappa_r(\widehat \T_{+m},d^{-m}\Pi^*[\Gamma_{-m}])$ tends to 0. This is a contradiction. The lemma follows.
\endproof

If $E_z$ and $E_w$ are two tangent subspaces at $z$ and $w$ in $\P^k$, of dimension respectively $k-p$ and $p$, then $E_z\times E_w$ is a tangent subspace of $\P^k\times\P^k$ at $(z,w)$, of dimension $k$. This construction induces a natural embedding of $\Gr(\P^k,k-p)\times \Gr(\P^k,p)$ into $\Gr(\P^k\times \P^k,k)$. For simplicity, we identify $\Gr(\P^k,k-p)\times \Gr(\P^k,p)$ with its image. Since $F$ is a product map, $\Gr(\P^k,k-p)\times \Gr(\P^k,p)$ can be seen as an invariant submanifold of $\Gr(\P^k\times \P^k,k)$.

\begin{lemma} \label{lemma_support_T+}
The current $\widehat \T_+$ is supported by $\Gr(\P^k,k-p)\times \Gr(\P^k,p)$.
\end{lemma}
\proof
By Proposition \ref{prop_F_laminar}, there is a measurable web associated with $\T_+$ which can be lifted to a measurable web of $\widehat \T_+$ in $\Gr(\P^k\times\P^k,k)$.
Since $\T_+=T_+\otimes T_-$, we have $\T_+\wedge \pi_1^*(\omega_\FS^{k-p+1})=0$ and $\T_+\wedge \pi_2^*(\omega_\FS^{p+1})=0$. Therefore, if $L$ is a generic lame of
$\T_+$, its tangent space at a regular point is the product of a tangent subspace of dimension $k-p$ of the first factor $\P^k$ and a tangent subspace of dimension $p$ of the second one. So 
the regular part of $L$ is locally a product of a manifold of dimension $k-p$ in $\P^k$ and another of dimension $p$. 
It follows that $\widehat L$ is contained in $\Gr(\P^k,k-p)\times \Gr(\P^k,p)$.
The lemma follows.
\endproof

We will see later that the currents $\widehat T_+$, $\widehat T_-$ in the following lemma are lifts of $T_+$, $T_-$ to $\Gr(\P^k,k-p)$ and $\Gr(\P^k,p)$ respectively, see Section \ref{section_woven} for the definition. 

\begin{lemma} \label{lemma_Ts}
Let $\widehat \T_s$ and $\widehat \T_u$ denote the push-forward of $\widehat \T_+$ to $\Gr(\P^k,k-p)\times\P^k$ and $\P^k\times \Gr(\P^k,p)$ respectively. Then there are positive closed currents $\widehat T_+$ on $\Gr(\P^k,k-p)$ and $\widehat T_-$ on $\Gr(\P^k,p)$ such that 
$\widehat \T_s=\widehat T_+\otimes T_-$ and $\widehat \T_u=T_+\otimes \widehat T_-$. 
\end{lemma}
\proof
We prove the existence of $\widehat T_+$. The case of $\widehat T_-$ can be obtained in the same way.
Let $\widehat\T_{+,n}$ be a limit value of $d^{-n_i+2n} [\widehat\Gamma_{n_i-2n}]$ when $i\to\infty$. Then $\widehat \T_{+,n}$ is also supported by $\Gr(\P^k,k-p)\times \Gr(\P^k,p)$ and $\widehat \T_+=d^{-2n} (\widehat F^n)^\bullet (\widehat \T_{+,n})$. 
Let $\widehat \T_{s,n}$ denote the push-forward of $\widehat\T_{+,n}$ to $\Gr(\P^k,k-p)\times\P^k$. 
Define $h:=(\widehat f,f^{-1})$ the product map on $\Gr(\P^k,k-p)\times\P^k$. We have $\widehat \T_s=d^{-2n} (h^n)^\bullet(\widehat \T_{s,n})$. Recall that $\Gr(\P^k,k-p)$ is bi-rational to $\P^k\times \G$. We still denote by $\widehat f$ the canonical lift of $f$ to 
$\P^k\times\G$. 

Recall that $d^{-n}[\Gamma_n]$ converges to $T_+\otimes T_-$  which is supported in $\overline K_+\times\overline K_-$. Therefore, 
the support of $\widehat \T_+$ is contained in the fibers over $\overline K_+\times\overline K_-$. 
We can apply Proposition  \ref{prop_product_map} to $\widehat\T_{s,n}$ and then Proposition 
\ref{prop_current_tensor}  (we have to permute the factors of $X\times Y$ in those propositions).
The current
$R$ defined as in Proposition \ref{prop_current_tensor} is supported by $\overline K_-$. Theorem \ref{th_unique_Green} applied to $f^{-1}$ implies that $R$ is necessarily a multiple of $T_-$.  Since $T_-$ is extremal, Proposition \ref{prop_current_tensor} implies the result.
\endproof

Denote by $\Pi_s$, $\Pi_u$, $\pi_s$ and $\pi_u$ the canonical projections from $\Gr(\P^k,k-p)\times\Gr(\P^k,p)$ onto $\Gr(\P^k,k-p)\times \P^k$, $\P^k\times \Gr(\P^k,p)$, $\Gr(\P^k,k-p)$ and $\Gr(\P^k,p)$ respectively. The following lemma shows that
$\widehat T_+$, $\widehat T_-$ are lifts of $T_+$, $T_-$ to $\Gr(\P^k,k-p)$ and $\Gr(\P^k,p)$ respectively.

\begin{lemma} \label{lemma_Ts_lim}
Let $\widehat \Omega_p$ be the standard lift of $\omega_\FS^p$ to $\Gr(\P^k,k-p)$ defined in Example \ref{ex_Fubini_bis}. Then $d^{-n_i} (\widehat f^{n_i})^*(\widehat\Omega_p)$  converges to 
$\widehat T_+$  and a similar result holds for 
$\widehat T_-$. 
\end{lemma}
\proof
Let $L$ be a generic projective subspace of dimension $k-p$ in  $\P^k$. 
By definition of  $\widehat\Omega_p$, 
we only have to prove that  $d^{-n_i}(\widehat f^{n_i})^\bullet [\widehat L]$ converges  to $\widehat T_+$. 
We will identify $\Gr(\P^k\times\P^k,k)$ with  $\P^k\times\P^k\times\GG$ via the natural bi-rational map and we still denote by $\Pi$ the projection onto $\P^k\times\P^k$. 

Fix a generic  $L$ as in Corollary \ref{cor_inter_woven} (we change the notation $H_\xi$ by $L$). So up to extracting a subsequence of $(n_i)$, we can assume that the intersections $d^{-n_i}[\widehat\Gamma_{n_i}]\wedge [\P^k\times L\times\GG]$ and $\widehat\T_+\wedge [\P^k\times L\times\GG]$ are well-defined and $d^{-n_i}[\widehat\Gamma_{n_i}]\wedge[\P^k\times L\times\GG]$ converges to $\widehat\T_+\wedge [\P^k\times L\times\GG]$.
Moreover, we can associate to $d^{-n_i}[\widehat\Gamma_{n_i}]$ and $\widehat\T_+$ the measurable webs which are described in that corollary. We also denote them by $\nu_{n_i}$, $\nu_{n_i}[s]$, $\nu$ and $\nu[s]$. 

\medskip\noindent
{\bf Claim.} Any limit value $S$ of  $d^{-n_i}(\widehat f^{n_i})^\bullet [\widehat L]$ is larger than or equal to  $\widehat T_+$.

\medskip

Assuming the claim, we first complete the proof of the lemma. 
We have seen in the proof of Theorem \ref{th_lam_Green} that $S$ is a lift of $T_+$ to $\Gr(\P^k,k-p)$
and the projection of a generic lame of $S$ to $\P^k$ is non-degenerate, i.e. locally of maximal dimension $k-p$. On the other hand, $\Pi_*(\widehat\T_+)$ is the limit of $d^{-n}[\Gamma_n]$ which is equal to $\T_+=T_+\otimes T_-$. We deduce that
the push-forward of $\widehat\T_s$ to $\P^k\times\P^k$ is also equal to $\T_+$ and hence
the push-forward of $\widehat T_+$  to $\P^k$ is equal to $T_+$. Since the same property holds for $S$, 
the current $S-\widehat T_+$, which is positive according to the claim, has horizontal dimension $<k-p$ with respect to the projection onto $\P^k$. But since generic lames of $S$ have horizontal dimension $k-p$, we necessarily have $S-\widehat T_+=0$. 
The lemma follows.
\hfill $\square$

\medskip\noindent
{\bf Proof of the claim.} The idea is to show that we can approach $\widehat T_+$ by lames in $\widehat f^{-n_i}(\widehat L)$ or more precisely by currents smaller than or equal to $d^{-n_i}(\widehat f^{n_i})^\bullet [\widehat L]$. 
Define $L_{n_i}$ as the closure in $\P^k$ of $f^{-n_i}(L)\cap \C^k$ and $L_{(n_i)}$ the closure in $\P^k\times \P^k$ of $\pi_2^{-1}(L)\cap\Gamma_{n_i}\cap \C^k\times\C^k$. The last analytic set is the family of points $(z,f^{n_i}(z))$ with $z\in f^{-n_i}(L)\cap\C^k$. So we have $\pi_1(L_{(n_i)})=L_{n_i}$. Define also $\widehat L_{n_i}$ as the lift of $L_{n_i}$ to $\Gr(\P^k,k-p)$ and $\L_{(n_i)}:=\Pi^{-1}(L_{(n_i)})\cap\widehat\Gamma_{n_i}$.  

The discussion before the claim on the intersection of currents with $\P^k\times L\times \GG$ implies that $d^{-n_i}[\L_{(n_i)}]\to \widehat \T_+\wedge [\P^k\times L\times\GG]$. The last intersection is supported by $\Gr(\P^k,k-p)\times\Gr(\P^k,p)$. 
By Lemma \ref{lemma_Ts}, 
its push-forward to $\Gr(\P^k,k-p)\times \P^k$ is equal to 
 $\widehat T_+\otimes ([L]\wedge T_-)$. Therefore, its push-forward to $\Gr(\P^k,k-p)$ is equal to the one of 
  $\widehat T_+\otimes ([L]\wedge T_-)$
 and hence equal to $\widehat T_+$ because  $[L]\wedge T_-$ is a probability measure. 
So in order to obtain the claim, it is enough to approach  $\widehat \T_+\wedge [\P^k\times L\times\GG]$ by woven currents on 
 $\Gr(\P^k,k-p)\times\Gr(\P^k,p)$ whose push-forwards to $\Gr(\P^k,k-p)$ are bounded by $d^{-n_i} [\widehat L_{n_i}]$. 

For this purpose,  we cannot directly use the currents $d^{-n_i}[\L_{(n_i)}]$ because they are not supported by $\Gr(\P^k,k-p)\times\Gr(\P^k,p)$ and they admit no natural push-forward to $\Gr(\P^k,k-p)$. Recall that we however have $d^{-n_i}[\L_{(n_i)}]\to \widehat \T_+\wedge [\P^k\times L\times\GG]$.
The idea now is to modify the lames of $d^{-n_i}[\L_{(n_i)}]$ without changing their limits in order to get convenient lames in $\Gr(\P^k,k-p)\times\Gr(\P^k,p)$. We have to make sure that the push-forwards of those lames to $\Gr(\P^k,k-p)$ satisfy the desired property.

Let $\widehat \Lambda$ be a graph corresponding to a generic point with respect to the measure $\nu[s]$ for some $s$. 
It is contained in $\Gr(\P^k,k-p)\times \Gr(\P^k,p)$.
Since $\widehat \T_+$ is a lift of $\T_+$, the projection $\Lambda:=\Pi(\widehat\Lambda)$ of $\widehat\Lambda$ to $\P^k\times\P^k$ corresponds to a lame of $\T_+$. 
We can refine the webs if necessary in order to get $\Lambda$ smooth and contained in $\C^k\times\C^k$.
Since $\T_+=T_+\otimes T_-$, we have seen that $\Lambda$ is locally a product $\Lambda_+\times\Lambda_-$ of a manifold $\Lambda_+$ of dimension $k-p$ in $\P^k$ with another $\Lambda_-$ of  dimension $p$. We can also write locally 
$\widehat\Lambda=\widehat\Lambda_+\times\widehat\Lambda_-$, where $\widehat\Lambda_+,\widehat\Lambda_-$ are the lifts of $\Lambda_+,\Lambda_-$ to $\Gr(\P^k,k-p)$ and $\Gr(\P^k,p)$ respectively. 

The projection $\widehat\Lambda_s:=\Pi_s(\widehat\Lambda)$ of $\widehat\Lambda$ to $\Gr(\P^k,k-p)\times\P^k$ is locally $\widehat\Lambda_+\times\Lambda_-$ and corresponds to a lame of $\widehat\T_s$. We see that the projection of $\widehat\Lambda\cap (\P^k\times L\times \GG)$ and the projection of $\widehat \Lambda_s\cap (\Gr(\P^k,k-p)\times L)$ to $\Gr(\P^k,k-p)$ are both equal to the lift of $\pi_1(\Lambda\cap \pi_2^{-1}(L))$ to $\Gr(\P^k,k-p)$. So each connected component $\lambda$ of $\widehat\Lambda\cap (\P^k\times L\times\GG)$ is the product of the lift of $\pi_1(\Pi(\lambda))$ to  $\Gr(\P^k,k-p)$ with a point in $\Gr(\P^k,p)$. 

Consider now a graph $\widehat\Lambda'$ corresponding to a generic point with respect to $\nu_{n_i}[s]$ which is close enough to $\widehat \Lambda$. It is an open subset of $\widehat\Gamma_{n_i}$.  Recall that the fat extension of $\widehat\Lambda'$  is also close to the fat extension of $\widehat\Lambda$. This insures that the approximation below is uniform on the  unextended lames. 
Consider a connected component $\widehat\lambda'$  of $\widehat\Lambda'\cap [\P^k\times L\times\GG]$. This is an open subset of $\L_{(n_i)}$.
We deduce from the above properties of $\widehat\Lambda$  that $\widehat\lambda'$ can be approximated by the product of the lift of $\pi_1(\Pi(\widehat\lambda'))$ to $\Gr(\P^k,k-p)$ with a point in $\Gr(\P^k,p)$. 
We need here the fact that $\P^k\times L\times\GG$ is transverse to the fat extension of $\widehat\Lambda$ which is guaranteed by Corollary \ref{cor_inter_woven} used just before the statement of the claim. 

Denote by $\widehat\lambda''$ the above product which is an approximation of $\widehat\lambda'$. The choice is not unique but the projection of $\widehat\lambda''$ in   $\Gr(\P^k,k-p)$ is always the lift of $\pi_1(\Pi(\widehat\lambda'))$.
Observe that since $\pi_1(\Pi(\widehat\lambda'))$ is an open subset of the variety $L_{n_i}$, the projection of $\widehat \lambda''$ to $\Gr(\P^k,k-p)$ is a lame in $\widehat L_{n_i}$. 
The approximation can be controlled uniformly on graphs $\widehat \Lambda'$ close enough to  $\widehat\Lambda$. Therefore, we can apply Lemma \ref{lemma_cv_inter} to $\nu_{n_i}[s]$ and $\nu[s]$ for each $s$. 
By considering  the projection onto $\Gr(\P^k,k-p)$, we deduce that every limit value of $d^{-n_i}[\widehat L_{n_i}]$ is larger than or equal to $\widehat T_+$. This is the claim.
\hfill $\square$

\medskip

\noindent
{\bf End of the proof of Proposition \ref{prop_a_weak_intersection}.}
Recall that when $\widehat \T_+\curlywedge \Pi^*[\Delta]$ is well-defined, $\widehat \T_+$ admits a unique tangent current along $\Pi^{-1}(\Delta)$ and the above intersection is the so-called shadow of this tangent current onto $\Pi^{-1}(\Delta)$, see  \cite{DS12}. Let $\widehat\nu$ denote the shadow of a tangent current of $\widehat \T_+$ along $\Pi^{-1}(\Delta)$. 
By Lemma \ref{lemma_kappa_T+}, this is a positive measure with support in $\Gr(\P^k,k-p)\times \Gr(\P^k,p)$.
It is enough to prove that $\widehat\nu=\widehat\mu^\Delta$. Indeed, in this case all tangential currents are necessarily vertical and  equal to the pull-back of $\widehat\mu^\Delta$. 

Define $\mu^\Delta:=\T_+\wedge [\Delta]$. 
Since $\T_+=T_+\otimes T_-$, it is not difficult to see that $\mu^\Delta=(\pi_{1|\Delta})^*(\mu)$, where $\mu:=T_+\wedge T_-$ is the Green measure of $f$. 
Since $\Pi_*(\widehat \T_+)=\T_+$, we deduce that $\Pi_*(\widehat\nu)=\mu^\Delta$. In particular, $\widehat\nu$ is a probability measure. By definition of $\widehat\mu^\Delta$, we also have  $\Pi_*(\widehat \mu^\Delta)=\mu^\Delta$. So we need to prove that for $\mu$-almost every $z\in\P^k$ the conditional measures 
$\langle \widehat \nu|\Pi|(z,z)\rangle$ and $\langle \widehat\mu^\Delta|\Pi|(z,z)\rangle$ of $\widehat\nu$ and $\widehat\mu^\Delta$ with respect to $\Pi$ are equal. Note that the second conditional measure is the Dirac mass at the point $(z,z,[E_s(z)\times E_u(z)])$ in the fiber $\{(z,z)\}\times\GG$ of $\Pi$. 

Since $\widehat\T_s=(\Pi_s)_*(\widehat \T_+)$, if $\Pi_+:\Gr(\P^k,k-p)\times\P^k\to \P^k\times\P^k$ denotes the canonical projection, then $(\Pi_s)_*(\widehat\nu)=\widehat\mu_s$, where $\widehat\mu_s:=\widehat\T_s\curlywedge\Pi_+^*[\Delta]$
provided that the last intersection exists. 
Observe that $\Pi_+^{-1}(\Delta)$ is the graph  of the canonical projection $\pi:\Gr(\P^k,k-p)\to\P^k$.
Since $\widehat \T_s=\widehat T_s\otimes T_-$ and $T_-$ is a power of a positive closed $(1,1)$-current with continuous potentials on $\C^k$, it is not difficult to see that the last intersection exists and is equal to  
$\widehat T_s\wedge \pi^*(T_-)$ if we identify $\Pi_+^{-1}(\Delta)$ with $\Gr(\P^k,k-p)$ in the canonical way, see Lemma 5.11 in \cite{DS12}.
By Lemma \ref{lemma_Ts_lim} and Proposition \ref{prop_main_inter} with $q=0$ and $\alpha=\widehat\Omega_p$, the last measure is equal to $\mu_+$.
The constant $c$ in this proposition is 1 since $\widehat\nu$ is a probability measure. 

We deduce that the conditional measure $\langle \widehat\mu_s |\Pi_+|(z,z)\rangle$ is equal to the Dirac mass at the point $((z,[E_s(z)]),z)$. It follows that the push-forward of $\langle \widehat \nu|\Pi|(z,z)\rangle$ to $\Gr(\P^k,k-p)\times\P^k$ is equal to the Dirac mass at the point 
$((z,[E_s(z)]),z)$. In the same way, we obtain that the push-forward of $\langle \widehat \nu|\Pi|(z,z)\rangle$ to $\P^k\times \Gr(\P^k,p)$ is equal to the Dirac mass at the point
$(z,(z,[E_u(z)]))$. We conclude that $\langle \widehat \nu|\Pi|(z,z)\rangle$ is equal to the Dirac mass at $(z,z,[E_s(z)\times E_u(z)])$.
This completes the proof of the proposition.
\hfill $\square$

\medskip

Recall that if $\Gamma$ is an analytic subset of dimension $k$ in $\P^k\times\P^k$ we define $\widetilde\Gamma$ 
as the closure of the set of points $(x,[v])$ in $\Gr(\P^k\times\P^k,k)$ with $x$ in the smooth part of $\Gamma$ and $v$ a non-zero complex tangent $k$-vector to $\P^k\times\P^k$ at $x$ which is not transverse to $\Gamma$. Let $\Sigma$ denote the set of points $(x,[v],[w])$, where $x\in\P^k\times\P^k$, $v$ and $w$ are non-zero complex tangent $k$-vectors of $\P^k\times\P^k$ at $x$. This is a smooth complex manifold. Over $\C^k\times\C^k$, we can identify it with $\C^k\times\C^k\times\GG\times \GG$. 
Denote by $\Sigma'$ the analytic subset of points $(x,[v],[w])\in\Sigma$ such that $v$ is not transverse to $w$. Over $\C^k\times\C^k$, it can be identified with the product of $\C^k\times\C^k$ with a hypersurface of $\GG\times \GG$, possibly with singularities.

Let $p_1$ and $p_2$ denote the projections from $\Sigma$ to $\Gr(\P^k\times\P^k,k)$ given by 
$p_1(x,[v],[w]):=(x,[v])$ and $p_2(x,[v],[w]):=(x,[w])$. So $p_1$ and $p_2$ define two fibrations on $\Sigma$ which are locally trivial with fibers isomorphic to $\GG$. They also define two fibrations on $\Sigma'$ which are locally trivial but the fibers may be singular. 
We have
$$\widetilde\Gamma=p_1\big(p_2^{-1}(\widehat\Gamma)\cap\Sigma'\big).$$
Define $\widetilde\T_+$ as the limit of $d^{-n_i}[\widetilde\Gamma_{n_i}]$. The last formula implies that this limit exists and is equal to 
$$\widetilde\T_+=(p_1)_*\big(p_2^*(\widehat\T_+)\wedge[\Sigma']\big).$$
The last wedge-product exists because the restriction of $p_2$ to $\Sigma'$ is a (possibly singular) fibration, i.e. locally a product. 

We now prove the crucial point in our approach for the main theorem. It says that $\Gamma_n$ is mostly transverse to $\Delta$ when $n\to\infty$. 

\begin{corollary} \label{cor_density_T}
The density between $\widetilde \T_+$ and $\widehat\Delta$ is zero.
\end{corollary}
\proof
By Proposition \ref{prop_a_weak_intersection}, the current $\widehat \T_+$ admits a unique tangent current along $\Pi^{-1}(\Delta)$ whose shadow on $\Pi^{-1}(\Delta)$ is $\widehat \mu^\Delta$. Recall that this tangent current can be defined as the limit of the images of $\widehat \T_+$ under local holomorphic  dilations in the normal directions to $\Pi^{-1}(\Delta)$. 
In the present situation, we can work over $\C^k\times\C^k$ where we can identify $\Sigma$ and $\Sigma'$ with products of $\C^k\times\C^k$ with some varieties. This allows us to use the obvious dilations given in a coordinate system on $\C^k\times\C^k$. 

So if $P:\Sigma\to\P^k\times\P^k$ is the canonical projection, 
we can also use the dilations in the normal directions to $P^{-1}(\Delta)$ induced by the obvious dilations in the normal directions of $\Delta$ in $\C^k\times\C^k$. Therefore, we deduce from
 Proposition \ref{prop_a_weak_intersection} that
$$(p_2^*(\widehat\T_+)\wedge[\Sigma']) \curlywedge [P^{-1}(\Delta)]=p_2^*(\widehat\mu^\Delta)\wedge [\Sigma'].$$
Using the above transform $\widehat\T_+\mapsto\widetilde \T_+$, we obtain
$$\widetilde\T_+\curlywedge[\Pi^{-1}(\Delta)]=(p_1)_*\big(p_2^*(\widehat\T_+)\wedge[\Sigma'] \big) \curlywedge [\Pi^{-1}(\Delta)]=(p_1)_*\big(p_2^*(\widehat\mu^\Delta)\wedge [\Sigma']\big)$$
on $\Gr(\P^k\times\P^k,k)$. 

The right hand side of the last identity is a positive closed current with support in $\Pi^{-1}(\Delta)$. We denote it by $R$.
In order to compute $R$, we first replace $\widehat\mu^\Delta$ with a Dirac mass $\delta$ at a point $(x,[v])$ in $\Gr(\P^k\times\P^k,k)$ with $x\in\C^k\times\C^k$. We identify $\Gr(\P^k\times\P^k,k)$ with $\C^k\times\C^k\times\GG$ over $\C^k\times\C^k$. 
If $\GG(v)$ denotes the hypersurface of $[w]\in\GG$ with $w$ not transverse to $v$, then  $(p_1)_*\big(p_2^*(\delta)\wedge [\Sigma']\big)$ is equal to the current of integration on $\{x\}\times \GG(v)$. 

For $\mu$-almost every $z$, we denote by $\GG(z)$ the hypersurface $\GG(v)$ with a vector $v$ defining  $E_s(z)\times E_u(z)$. The last vector space is identified with a subspace of the tangent space of $\P^k\times\P^k$ at $(z,z)$. 
Using the definition of $\widehat\mu^\Delta$, 
we deduce that the above current $R$ is equal to
$$R=\int [\{(z,z)\}\times \GG(z)] d\mu(z).$$
The tangent current to $\widetilde\T_+$ along $\Pi^{-1}(\Delta)$ is the pull-back of $R$ to the normal vector bundle of $\Pi^{-1}(\Delta)$ in $\Gr(\P^k\times\P^k,k)$. We denote it by $S$. 

Observe that since $E_s(z)\cap E_u(z)=\{0\}$, the two vector spaces $E_s(z)\times E_u(z)$ and $\Delta\cap \C^k\times\C^k$ only intersect at the point $(z,z)$. So
$E_s(z)\times E_u(z)$ is transverse to $\Delta$ and  
each manifold $\{(z,z)\}\times \GG(z)$ is disjoint from $\widehat\Delta$.  Hence, the tangent current of $S$ along $\widehat\Delta$ is  zero. Therefore, Proposition 4.13
 in \cite{DS12} implies that the tangent current of 
$\widetilde\T_+$ along $\widehat\Delta$ is zero. This completes the proof of the corollary.
\endproof

In what follows, denote by $(x,y)$ the standard coordinates in $\C^k\times\C^k$. Define also $z:=x-y$ and $w:=x$. So the diagonal $\Delta$ is given by $x=y$ or by $z=0$. 

\begin{lemma} \label{lemma_near_Delta}
Let $R_0>0$ be a fixed constant. Then
if $R>0$ is a constant large enough then $\Gamma_n\cap\{\|z\|\leq R_0\}$ is contained in $\{\|w\|< R\}$ for every $n\geq 1$. 
\end{lemma}
\proof
Fix a neighbourhood $U_+$ of $I_+$ and a neighbourhood $U_-$ of $I_-$ small enough such  that $\overline U_+\cap \overline U_-=\varnothing$. Since $I_+$ is attractive for $f^{-1}$, we can choose $U_+$ so that $f(\P^k\setminus U_+)\subset \P^k\setminus U_+$. We can also assume that $f^{-1}(\P^k\setminus U_-)\subset \P^k\setminus U_-$. 
Fix  a constant $R>0$ large enough. Consider a point $(x,y)$ such that $\|z\|\leq R_0$ and $\|w\|\geq R$. We have $\|x\|\geq R-R_0$, $\|y\|\geq R-R_0$ and $\|x-y\|\leq R_0$. Therefore, either $x\not\in U_+$ or $y\not\in U_-$. We have to show that $(x,y)\not\in\Gamma_n$. Assume without loss of generality that $x\not\in U_+$. It suffices  to prove that $y\not=f^n(x)$. 

Let $f_+$ denote the homogeneous part of maximal degree $d_+$ of $f$. Then the closure of $\{f_+=0\}$ is an analytic subset of $\P^k$ whose intersection with the hyperplane at infinity is equal to $I_+$. For $x$ out of a neighbourhood of this analytic set, we have that $\|f(x)\|\geq c\|x\|^{d_+}$ for some constant $c>0$. In particular, this holds for $x$ as above when the neighbourhood of $\{f_+=0\}$ is small enough. So we have $\|f(x)\|>2\|x\|$. This together with the inclusion $f(\P^k\setminus U_+)\subset \P^k\setminus U_+$ imply by induction that $\|f^n(x)\|\geq 2^n\|x\|>\|y\|$. The lemma follows.
\endproof

For simplicity, using a linear change of coordinates, we can assume that $R=1$ satisfies the above lemma for some constant $R_0>0$. Define $z':=\lambda z$ and $w':=w$ for some large constant $\lambda>4/R_0$ that will be made precise later. We will apply Proposition \ref{prop_branch_bis} to the restriction of $\Gamma_n$ in the domain $\U:=4\B_k\times 3\B_k$ with respect to the coordinates $(z',w')$. 

\begin{lemma} \label{lemma_degree_Gamma}
Let $d_n$ denote the number of periodic points of period $n$ of $f$ in $\C^k$ counted with multiplicity. Then
the restriction of $\Gamma_n$ to  $\U$ is contained in $4\B_k\times\B_k$ and is a ramified covering of degree $d_n$
over the  factor $4\B_k$. Moreover, we have $d_n=d^n+o(d^n)$.  
\end{lemma}
\proof
It is clear that $\Gamma_n\cap \U\subset 4\B_k\times \B_k$ so $\Gamma_n\cap\U$ is a ramified covering over $4\B_k$. 
Its degree is equal to the number of points in the intersection $\Gamma_n\cap\Delta\cap (\C^k\times\C^k)$ counted with multiplicity. So this degree is equal to $d_n$.

Let $\pi:\C^k\times\C^k\to\C^k$ denote the projection $(x,y)\mapsto x-y$. If $\nu$ is a smooth probability measure with compact support in a small neighbourhood of $0$ in $\C^k$ then $d_n$ is the mass of the measure $[\Gamma_n]\wedge \pi^*(\nu)$.  
Since $d^{-n}[\Gamma_n]\to \T_+$, the sequence $d^{-n}d_n$ converges to the mass of $\T_+\wedge \pi^*(\nu)$. 
In particular, this mass does not depend on the choice of $\nu$. In order to compute this mass, we take a sequence of measures $\nu_n$ converging to the Dirac mass at 0. Since
$\T_+=T_+\otimes T_-$ 
and $T_\pm$ are powers of positive closed $(1,1)$-currents with continuous potentials,  $\T_+\wedge \pi^*(\nu_n)$ converge
to $(T_+\otimes T_-)\wedge [\Delta]$ which is the probability measure $\mu^\Delta$. We conclude that $d^{-n}d_n\to 1$. This completes the proof of the lemma.
\endproof

Recall that we are using the coordinate systems $z':=\lambda z=\lambda(x-y)$ and $w':=w:=x$ in $\C^k\times\C^k$. 

\begin{lemma} \label{lemma_branch_Gamma}
Let $0<\delta<1$ be a fixed constant. Fix also a constant $\lambda>0$ large enough depending on $\delta$. Then for $n$ large enough, $\Gamma_n\cap (\B_k\times 3\B_k)$ admits at least $(1-\delta^2)d^n$ connected components which are graphs over $\B_k$.
\end{lemma}
\proof
We want to apply Proposition \ref{prop_branch_bis}. Fix an $n$ large enough and  define $\Gamma:=\Gamma_n\cap \U$. If $\widetilde\Gamma^\star$ is defined as in  Proposition \ref{prop_branch_bis}, we have to show that $d^{-n}\|\widetilde\Gamma^\star\|$ is as small as we want when $\lambda$ is large enough.

We first consider $\Gamma_n$ and $\Delta$ in the coordinate system $(z,w)$ and recall that $\Gamma_n\cap \{\|z\|\leq R_0\}$ is contained in $\{\|w\|\leq 1\}$. Define $\U_0:=\{\|z\|<4, \|w\|<3\}$ and $\Delta_0:=\Delta\cap \U_0$. 
Recall also that $\Pi:\Gr(\P^k\times\P^k,k)\to\P^k\times\P^k$ is the canonical projection and $\Pi^{-1}(\U_0)$ is identified with the product $\U_0\times\GG$, where $\GG$ is the Grassmannian of linear subspaces of dimension $k$ in $\C^k\times\C^k$. Consider the square matrix $A$ with complex coefficients of size $k\times k$. The linear subspace $z=Aw$ corresponds to a point in $\GG$. So we can use $A$ for affine coordinates of a Zariski open subset $\GG_0$ of $\GG$.  In these coordinates, the lift $\widehat \Delta_0$ of $\Delta_0$ to $\Gr(\P^k\times\P^k,k)$ is identified with $\Delta_0\times\{0\}$. 

The coordinate change $(z',w')=(\lambda z,w)$ on $\U_0$ induces the coordinate change $(z',w',A'):=(\lambda z, w, \lambda A)$ on $\U_0\times\GG_0$ and can be seen as the dilation along the normal directions to $\widehat\Delta_0$, as we have used when we defined the tangent currents in \cite{DS12}. By Corollary \ref{cor_density_T}, the density between $\widetilde \T_+$ and $\widehat\Delta$ is zero. Thus, since $\lambda$ is large enough,  the mass of $\widetilde\T_+$ in $\W:=\{\|z'\|< 4, \|w'\|< 3,\|A\|< 1\}$ with respect to the standard Euclidean metric associated with these coordinates, is as small as we want.
By definition of $\widehat\T_+$ and $\widetilde\T_+$, we have $d^{-n_i}[\widetilde\Gamma_{n_i}]\to \widetilde\T_+$. 
Therefore, the mass of $d^{-n_i}[\widetilde\Gamma_{n_i}]$ in $\W$ is as small as we want when $i$ is large enough. The property holds for every choice of $\widehat\T_+$. Since $n$ is large enough, we deduce that for $\Gamma,\widetilde\Gamma^\star$ defined above, $d^{-n}\|\widetilde\Gamma^\star\|$ is as small as we want. This completes the proof of the lemma.
\endproof

Recall that $\pi_1,\pi_2$ are the canonical projections from $\P^k\times\P^k$ onto its factors. 
Let $\Gamma^{(j)}_n$ denote one of the graphs obtained in Lemma \ref{lemma_branch_Gamma}. 
It is the graph of $f^n$ over the domain $\pi_1(\Gamma^{(j)}_n)$.
It intersects $\Delta$ at a unique point $(a^{(j)}_n,a^{(j)}_n)$, where $a^{(j)}_n$ is a periodic point of period $n$. The eigenvalues of the differential $Df^n(a_n^{(j)})$ of $f^n$ at $a^{(j)}_n$ do not depend on the local coordinate system at $a_n^{(j)}$. Fix a constant $0<\epsilon<1$ as in the introduction.

\begin{proposition} \label{prop_branch_Gamma}
Let $\delta$ and $\lambda$ be as in Lemma \ref{lemma_branch_Gamma}. Then for $n$ large enough there are at least 
$(1-\delta)d^n$ graphs $\Gamma_n^{(j)}$ such that $Df^n(a_n^{(j)})$ admits exactly $p$ eigenvalues with modulus $\geq (d_+-\epsilon)^{n/2}$ and $k-p$ eigenvalues with modulus $\leq (d_--\epsilon)^{-n/2}$, counted with multiplicity.
\end{proposition}
\proof
Using the action of $f^n$ on the cohomology of $\P^k$, we have 
 for $1\leq q\leq p$
$$d^{-n}\int_{\Gamma_n} \pi_1^*(\omega_\FS^{k-p+q})\wedge \pi_2^*(\omega_\FS^{p-q})=d^{-n}\int_{\P^k} \omega_\FS^{k-p+q}\wedge (f^n)^*(\omega_\FS^{p-q}) \leq d_+^{-n}$$
and for $1\leq q\leq k-p$
$$d^{-n}\int_{\Gamma_n} \pi_1^*(\omega_\FS^{k-p-q})\wedge \pi_2^*(\omega_\FS^{p+q})=
d^{-n}\int_{\P^k} (f^n)_*(\omega_\FS^{k-p-q})\wedge \omega_\FS^{p+q}\leq d_-^{-n}.$$

Observe that the graphs $\Gamma^{(j)}_n$ are contained in a fixed compact subset of $\C^k\times\C^k$. Moreover, on any fixed compact subset of $\C^k$, the standard K\"ahler form $i\ddbar \|x\|^2$  is comparable with $\omega_\FS$. Therefore, there are at least $(1-\delta)d^n$ graphs $\Gamma_n^{(j)}$ such that for $1\leq q\leq p$
$$\int_{\Gamma_n^{(j)}} (i\ddbar \|x\|^2)^{k-p+q}\wedge (i\ddbar\|y\|^2)^{p-q} \leq c d_+^{-n}$$
and for $1\leq q\leq k-p$
$$\int_{\Gamma_n^{(j)}} (i\ddbar\|x\|^2)^{k-p-q}\wedge (i\ddbar\|y\|^2)^{p+q} \leq c d_-^{-n},$$
where $c>0$ is a constant depending on $\delta,\lambda$ but independent of $n$. In what follows, we only consider the graphs satisfying these estimates.

It is convenient now to work with the coordinates $(z,w)=(x-y,x)$. Since $(z,w)=(\lambda^{-1}z',w')$,  $\Gamma_n^{(j)}$ is a graph in $\lambda^{-1}\B_k\times \B_k$ over $\lambda^{-1}\B_k$. Denote by $h:\lambda^{-1}\B_k\to\Gamma_n^{(j)}$ the canonical map. 

\medskip
\noindent
{\bf Claim.} There is a constant $c_1>0$ depending on $\delta,\lambda$ such that for $1\leq q\leq p$
$$\big\|h^*\big[(i\ddbar \|x\|^2)^{k-p+q}\big]\wedge h^*\big[(i\ddbar\|y\|^2)^{p-q}\big]\big\|_0\leq c_1d_+^{-n}$$
and for $1\leq q\leq k-p$
$$\big\|h^*\big[ (i\ddbar \|x\|^2)^{k-p-q}\big]\wedge h^*\big[(i\ddbar\|y\|^2)^{p+q}\big] \big\|_0\leq c_1d_-^{-n},$$
where $\|\cdot\|_a$ denotes the norm of a cotangent vector of maximal bi-degree $(k,k)$ at the point $a$. 

\medskip

Assuming the claim, we first complete the proof of the proposition. 
Define for simplicity $D:=Df^n(a^{(j)}_n)$.
Let $\gamma_1,\ldots,\gamma_k$ denote the eigenvalues of $D$ ordered so that $|\gamma_1|\leq \cdots\leq |\gamma_k|$. 
Let $l$ be an integer such that $|\gamma_i|\leq 1$ for $i\leq k-l$ and $|\gamma_i|\geq 1$ for $i>k-l$. We have either $l\leq p$ or $k-l\leq k-p$. Replacing $f$ with $f^{-1}$ if necessary, we can assume that $l\leq p$. 

The operator $D$ induces a linear automorphism, denoted by  $D^*$, of the cotangent space of $\C^k$ at $a_n^{(j)}$ and its exterior powers.  
Let $v$ be a unitary positive cotangent vector of bi-degree $(p-q,p-q)$
 which is an eigenvector associated with the eigenvalue $|\gamma_{k-p+q+1}|^{2}\ldots |\gamma_k|^{2}$ for some $1\leq q\leq p$.   Since $i\ddbar\|y\|^2$ is strictly positive, 
 $v$ is bounded by a constant times $i\ddbar\|y\|^2$.
 We deduce from the claim that 
 $$\big\|h^*\big[(i\ddbar\|x\|^2)^{k-p+q}\big]\wedge (\pi_2\circ h)^*(v)\big\|_0\leq c_2d_+^{-n}$$ 
 for some constant $c_2>0$. 
 
 Consider also the push-forward operator $(\pi_1\circ h)_*$ on forms of maximal bi-degree $(k,k)$ at 0. 
 This is the multiplication by the real Jacobian of the differential of $(\pi_1\circ h)^{-1}$ at $a_n^{(j)}$. This Jacobian is equal to $|\gamma_1-1|^{2}\ldots |\gamma_k-1|^{2}$ because $(\pi_1\circ h)^{-1}(x)=f^n(x)-x$. Consider the value of this operator acting at the vector inside the sign $\|\ \|_0$ of  the last inequality. We get 
 $$\big\|(i\ddbar\|x\|^2)^{k-p+q}\wedge D^*(v)\big\|_{a_n^{(j)}}\leq c_2d_+^{-n}|\gamma_1-1|^{2}\ldots |\gamma_k-1|^{2}.$$
 Hence 
 $$|\gamma_{k-p+q+1}|^2\ldots |\gamma_k|^2 \big\|(i\ddbar\|x\|^2)^{k-p+q}\wedge v\big\|_{a_n^{(j)}}\leq c_2d_+^{-n}|\gamma_1-1|^{2}\ldots |\gamma_k-1|^{2}.$$
 Since $v$ is unitary, we deduce that 
 $$|\gamma_{k-p+q+1}|^2\ldots |\gamma_k|^2 \leq c_3d_+^{-n}|\gamma_1-1|^{2}\ldots |\gamma_k-1|^{2}$$
 for some constant $c_3>0$.
 
 Consider first the case $l<p$ and take $q=p-l$. The choice of $l$ implies 
 that $|\gamma_i-1|\leq 2$ for $i\leq k-l$ and $|\gamma_i-1|\leq 2|\gamma_i|$ for $i\geq k-l+1$. This contradicts  the inequality obtained above when $n$ is large enough. So we have $l=p$. By taking $q=1$, we deduce from the same inequality that $|\gamma_{k-p+1}|\geq c_4d_+^{n/2}$ for some constant $c_4>0$. Finally, since $l=p$, we can apply the same arguments to $f^{-1}$ instead of $f$ and obtain that 
 $|\gamma_{k-p}^{-1}|\geq c_4d_-^{n/2}$ since the eigenvalues of $Df^{-1}(a_n^{(j)})$ are $\gamma_1^{-1},\ldots,\gamma_k^{-1}$.  
The proposition follows.
\endproof

\noindent
{\bf Proof of the claim.} We prove the first estimate. The second one is obtained in the same way. We have 
$$h^*\big[(i\ddbar \|x\|^2)^{k-p+q}\big]\wedge h^*\big[(i\ddbar\|y\|^2)^{p-q}\big]= h^*\big[(i\ddbar \|x\|^2)^{k-p+q}
\wedge (i\ddbar\|y\|^2)^{p-q}\big].$$
Denote by $\Theta$ this form. Since it is positive of maximal bi-degree, we can write 
$$\Theta(z)=\varphi(z) (idz_1\wedge d\overline z_1)\wedge \ldots\wedge (idz_k\wedge d\overline z_k),$$ 
where $\varphi(z)$ is a positive function. 
We have to show that $\varphi(0)\leq c_1d_+^{-n}$ for some constant $c_1>0$. 

The estimate given just before the claim implies that the integral of $\Theta$ on $\lambda^{-1}\B_k$ is smaller than $cd_+^{-n}$. Therefore, it is enough to check that $\varphi$ is a psh function. 
Observe that $\Theta$ is a finite sum of forms of type
$$(i\ddbar|g_1|^2)\wedge \ldots\wedge (i\ddbar|g_k|^2),$$
where $g_1,\ldots,g_k$ are holomorphic functions. The last form is equal to 
$$|\Jac(g_1,\ldots,g_k)|^2(idz_1\wedge d\overline z_1)\wedge \ldots\wedge (idz_k\wedge d\overline z_k),$$
where $\Jac$ denotes the complex Jacobian of a holomorphic map.
It is now clear that $\varphi$ is psh. The claim follows.
\hfill $\square$

\medskip\noindent
{\bf End of the proof of Theorem \ref{th_main}.}
Define 
$$\mu_n:=d^{-n}\sum_{a\in Q_n} \delta_a \quad \mbox{and} \quad \mu_n^\Delta:=d^{-n}\sum_{a\in Q_n} \delta_{(a,a)},$$
where $Q_n$ is as in the statement of the theorem and $\delta_{(a,a)}$ denotes the Dirac mass at the point $(a,a)$ in $\Delta$. By Lemma \ref{lemma_degree_Gamma}, the mass of $\mu_n$ is bounded and any limit value of $\mu_n$ is of mass at most equal to 1. So in order to prove the theorem, it is enough to consider the case of smallest sets $Q_n$, i.e. $Q_n= SP_n^\epsilon$. 

If $(n_i)$ is an increasing sequence of integers such that $\mu_{n_i}$ converges to a measure $\mu'$, we only have to show that $\mu'\geq \mu$. 
It is more convenient to work on $\Delta$. Denote by $\mu'^\Delta$ the limit of $\mu^\Delta_{n_i}$. We have to check that $\mu'^\Delta\geq\mu^\Delta$. 
For this purpose, it suffices to construct a positive measure $\mu''^\Delta$ such that $\mu''^\Delta\leq\mu^\Delta$, $\mu''^\Delta\leq \mu'^\Delta$ and $\|\mu''^\Delta\|\geq 1-\delta$ for every fixed constant $0<\delta<1$.  

Let $\lambda>0$ be the constant satisfying Proposition \ref{prop_branch_Gamma}. Define 
$$\S_{n_i}:=d^{-n_i}\sum_j [\Gamma_{n_i}^{(j)}].$$
By extracting a subsequence, we can assume that $\S_{n_i}$ converges to a current $\S$. 
Since $\Gamma_{n_i}^{(j)}$ are graphs of bounded holomorphic maps, we have
$\S_{n_i}\wedge [\Delta]$ and $\S\wedge [\Delta]$ are well-defined and $\S_{n_i}\wedge [\Delta]\to \S\wedge [\Delta]$.
The last convergence is guaranteed by the fact that if a sequence of bounded graphs converges in the sense of currents then it also converges locally uniformly. 

Define $\mu''^\Delta:=\S\wedge[\Delta]$. Since $d^{-n}[\Gamma_n]\to \T_+$, we have $\S\leq \T_+$. It follows that $\mu''^\Delta\leq\mu^\Delta$. 
By Proposition \ref{prop_branch_Gamma}, $\S_{n_i}\wedge [\Delta]$ is a positive measure smaller than or equal to $\mu_{n_i}^\Delta$. 
Thus, $\mu''^\Delta\leq \mu'^\Delta$. Proposition  \ref{prop_branch_Gamma} also implies that the mass of $\mu''^\Delta$ is at least equal to $1-\delta$. This completes the proof of the theorem.
\hfill $\square$

\begin{remark} \rm
The current $\S$ is constituted by lames of $\T_+=T_+\otimes T_-$. So its lames are locally a product of two manifolds. This property induces a product structure of $\mu$ by stable and unstable manifolds. 
\end{remark}

\begin{remark}\rm \label{rk_final}
Let $L_+$ and $L_-$ be analytic subsets of $\P^k$ of pure dimension $k-p$ and $p$ respectively. Assume that $L_+\cap I_-=\varnothing$ and $L_-\cap I_+=\varnothing$. Using the same method we can show that the points in $f^{-n}(L_+)\cap f^n(L_-)$ are equidistributed with respect to $\mu$ as $n\to\infty$, see also \cite{DS6,Dujardin}. For the proof, we replace the graphs $\Gamma_n$ with $f^{-n}(L_+)\times f^n(L_-)=F^{-n}(L_+\times L_-)$. The problem is simpler because the last analytic sets are already products of varieties and their lifts to $\Gr(\P^k\times\P^k,k)$ are contained in $\Gr(\P^k,k-p)\times\Gr(\P^k,p)$. 
\end{remark}

T.-C. Dinh, UPMC Univ Paris 06, UMR 7586, Institut de
Math{\'e}matiques de Jussieu, F-75005 Paris, France.\\
DMA, UMR 8553, Ecole Normale Sup\'erieure,
45 rue d'Ulm, 75005 Paris, France.\\
 {\tt
  dinh@math.jussieu.fr}, {\tt http://www.math.jussieu.fr/$\sim$dinh}

\

\noindent
N. Sibony,
Universit{\'e} Paris-Sud, Math{\'e}matique - B{\^a}timent 425, 91405
Orsay, France. {\tt nessim.sibony@math.u-psud.fr}

\end{document}